\documentclass[10pt,twocolumn,twoside]{IEEEtran} 
\usepackage{graphics}
\usepackage{graphicx}
\usepackage{epsfig,amssymb}
\usepackage{amsmath}
\usepackage{geometry}
\usepackage{mathrsfs}
\usepackage{color}
\usepackage{array}
\usepackage{amsfonts,booktabs}
\usepackage{lipsum,amsmath,multicol}
\usepackage{graphicx}
\usepackage{subfigure}
\usepackage{times}
\usepackage{pdfsync}
\usepackage{pgfplots}
\usepackage{pgfplotstable}
\usepackage[subnum]{cases}
\usepackage{filecontents}
\usepackage{multirow}
\usepackage{algorithm}
\usepackage{algpseudocode}
\usepackage{float}
\usepackage{setspace}
\newtheorem{theorem}{Theorem}
\newtheorem{lemma}{Lemma}
\newtheorem{definition}{Definition}
\newtheorem{corollary}{Corollary}

\newcommand{\beq}{\begin{equation}}
\newcommand{\eeq}{\end{equation}}
\newcommand{\beqa}{\begin{eqnarray}}
\newcommand{\eeqa}{\end{eqnarray}}

\newcommand{\paren}[1]{\left(#1\right)}

\newcommand{\field}[1]{\ensuremath{\mathbb{#1}}}
\newcommand{\abs}[1]{\left|#1\right|} % absolute value
\newcommand{\N}{\ensuremath{\field{N}}} % natural numbers
\newcommand{\R}{\ensuremath{\field{R}}} % real numbers
\newcommand{\Z}{\ensuremath{\field{Z}}} % integers
\newcommand{\1}{\ensuremath{\mathbf{1}}} % vector of all 1's
\newcommand{\I}[1]{\ensuremath{\mathsf{1}_{\left\{#1\right\}}}} % indicator function
\newcommand{\IC}[1]{\ensuremath{\mathsf{I}_{\left\{#1\right\}}}} % Constant indicator
\newcommand{\PR}[1]{\ensuremath{\mathsf{Pr}\left\{#1\right\}}} % probability with braces
\newcommand{\ES}[1]{\ensuremath{\mathsf{E}\left[#1 \right]}} %Expectation with square parentheses
\newcommand{\ENT}[1]{\ensuremath{\mathsf{h}\left[#1 \right]}} %Entropy
\newcommand{\e}[1]{\ensuremath{{\rm e}^{#1}}} %Exponents of e

\renewcommand{\vec}[1]{\ensuremath{\boldsymbol{#1}}}

\newcommand{\logp}[1]{\ensuremath{\log\paren{#1}}}

\newcommand{\norm}[2]{\left\Vert#1\right\Vert_{#2}}

\newcommand{\quant}[3]{\ensuremath{Q^{#2}_{#3}}}%quantizer
\newcommand{\CEP}[2]{\ensuremath{\mathsf{N}\left[\left.#1\right|#2\right]}}%conditional differential entropy power
\newcommand{\CENT}[2]{\ensuremath{\mathsf{h}\left[\left.#1\right|#2\right]}}%conditional differential entropy 
\newcommand{\CES}[2]{\ensuremath{\mathsf{E}\left[\left.#1\right|#2\right]}}%conditional expectation
\newcommand{\QM}[2]{\ensuremath{\hat{T}\left(#1,#2\right)}}%quantized map
\begin{document}
\title{Convergence Analysis of Quantized Primal-dual Algorithms in Network Utility Maximization Problems}
\author{Ehsan Nekouei, Tansu Alpcan, Girish N.~Nair, Robin J.~Evans 
\thanks{Department of Electrical and Electronic Engineering, The University of Melbourne, VIC 3010, Australia. E-mails: \{ehsan.nekouei,tansu.alpcan,gnair,robinje\}@unimelb.edu.au }
}
\maketitle
\thispagestyle{empty}
\begin{abstract}
This paper investigates the asymptotic and non-asymptotic behavior of the quantized primal-dual (PD) algorithm in network utility maximization (NUM) problems,
%%GN - To this end, we consider a NUM problems
in which a group of agents  maximize the sum of their individual concave objective functions under linear constraints. 
%%In each iteration of the algorithm, the primal variables are updated by the agents, whereas the dual variables are updated by a central network manager (NM) which has access to the parameters characterizing network-wide constraints. Also
%%At each time step, the agents and network manager have only access to the quantized versions of dual and primal variables, respectively.
 In the asymptotic scenario, we 
%%consider a class of quantization schemes called optimum achieving quantization scheme which allow the primal and dual variables converge to their optimal solution regardless of quantized communication between agents and the NM. 
use the information-theoretic notion of {\em differential entropy power} to establish universal bounds on the maximum exponential convergence rates of joint PD, primal and dual variables under optimum-achieving quantization schemes. These results provide trade-offs between the speed of exponential convergence, the agents' objective functions, the communication bit rates, and the number of agents and constraints. In the non-asymptotic scenario, we obtain lower bounds on the mean square distance of joint PD, primal and dual variables from the optimal solution at any time instant. These bounds hold regardless of the quantization scheme used.
%%imply that at a given time,  the PD algorithm cannot be arbitrarily close to the  optimal solutions.
%% and provide bounds on the closeness of the PD algorithm to the optimal solution for any finite time instance.
\end{abstract}

\section{Introduction}
\subsection{Motivation}
With continuing advances in  networking technology, our societies have become increasingly dependent on network-based technologies for performing everyday tasks. 
For example, consider data transfer using the internet, environmental monitoring  using wireless sensor networks, and online storage or computation in the ``cloud''. 
In all these applications, a limited number of resources, \emph{e.g.}, bandwidth, memory and CPU time, are shared among a group of networked devices, hereafter called \emph{agents}, 
to deliver the required service. As the quality of the delivered task is highly dependent on how the network resources are shared among the agents, 
resource allocation algorithms have become  vital components of these technologies. 

%As the number of agents in the network grows, it becomes impractical to optimize resource sharing  in a centralized way, for a number of reasons. 
%Firstly, solving the resource allocation problem becomes exceedingly hard as the number of agents becomes large. 
%GN - ORDER OF COMPLEXITY IN N? NP HARD?
%Secondly, for large and heterogeneous networks such as the internet, centralized optimization algorithms are difficult to implement,
 %due to the lack of complete information about the network. Thirdly, networks with centralized structure are more prone to cyber-attacks and failure. 

In the seminal work \cite{KMT98}, Kelly \emph{et al.} introduced the network utility maximization (NUM) approach, which provides decentralized frameworks, \emph{e.g.,} primal, dual and primal-dual (PD) decomposition methods, for solving large-scale resource allocation problems.    %In the last decade, substantial research effort has been dedicated to the design and analysis of decentralized optimization algorithms for solving  network utility maximization (NUM) problems. These research efforts resulted in a suite of decomposition techniques, \emph{e.g.,} primal-dual (PD), primal, dual decomposition methods, which allows large scale optimization problems, the rate allocation problem in capacity limited communication networks, to be solved in a decentralized way.  resource allocation problem since many practical resource sharing problems can be formulated as a NUM problem \emph{e.g.,} . since it allows many large scale resource sharing problems, Many resource allocation problems fall under the umbrella of network utility maximization (NUM) problems {social welfare} maximization problem in which a central entity called the \emph{network manager} (NM) seeks to maximize the sum of  objective functions of a group of agents sharing limited resources.   When the number of agents is small, the NM can efficiently solve the resource allocation problem using efficient optimization algorithms. However, in a large scale network with heterogeneous agents, \emph{e.g.,} internet,  centralized optimization algorithms are not applicable due to the lack of required information \cite{AB05}. 
In each decomposition method, the computational burden of solving the resource allocation problem is distributed among agents, and the task of information transfer between different agents is handled by an underlying communication network. The problem of devising efficient decomposition methods for NUM problems has been extensively studied in the past decade, \emph{e.g.,} see \cite{SS07} and references therein. 
Our aim in this paper is is to analyze the impact of quantized communications in NUM problems, using information-theoretic ideas.

\subsection{Related Work}
Although the performance of distributed optimization algorithms, and in particular NUM algorithms, under perfect communication networks is well understood, 
the investigation of the impact of imperfect communications on these optimization algorithms is relatively a new research area that has attracted much interests in recent years, \emph{e.g.,} see \cite{Nedic08}-\cite{YH14}.

Nedi\'{c} {\em et al.} \cite{Nedic08} considered a convex optimization problem, in which a set of agents collaboratively minimize a sum of individual objective functions. They proposed an averaging-based algorithm and studied its convergence rate under an infinite-level, uniform quantization scheme.  In \cite{Rabbat05}, the authors proposed an incremental algorithm for solving a convex optimization problem. They analyzed the convergence of the proposed algorithm under a uniform quantization scheme.% In each iteration of an incremental algorithm, only one agent will  update the optimization variables based on its local information, and then, it passes the updated values of the optimization variables to the next agent.

Yuan \emph{et al.} \cite{YXZR12} considered a constrained optimization problem in which a group of agents cooperate to minimize the sum of their local convex objective functions subject to a set of global constraints. They proposed a dual averaging algorithm and analyzed its convergence under uniform deterministic/stochastic quantization schemes. The authors in \cite{YH14} proposed a distributed sub-gradient algorithm for solving  an unconstrained multi-agent convex optimization problem, and studied its convergence under uniform zoom-in quantization. Finally, the authors of \cite{CL10} studied the problem of minimizing an upper bound on the distortion due to quantization in distributed iterative algorithms. They established  the optimality of different quantization structures under various distortion measures. Different from the literature discussed above, in this paper we study the speed of \emph{exponential} convergence of quantized PD algorithm in solving NUM problems. Moreover, our main results are independent of the structure of the underlying quantization scheme and hence can be applied to a more general class  than uniform quantizers.%GN-IN ALL THE WORK ABOVE DO WE NEED TO EXPLAIN CONVERGENCE RATES ACHIEVED AND GRAPH STRUCTURE ASSUMED? 

In \cite{NNA15-CDC}, we studied the convergence behavior of the PD algorithm in a quadratic NUM problem under quantized communications. %The convergence result in \cite{NNA15-CDC} is only valid for joint PD variables and quadratic objective functions.
 In the current paper, the objective functions of agents belong to the class of concave and twice continuously differentiable functions. This complicates our analysis as the PD update rule becomes non-linear in primal variables. Here, we study the impact of quantized communications on the convergence behavior of \emph{joint PD, primal and dual} variables in both \emph{asymptotic and non-asymptotic} regimes.% We note that the convergence analysis under general non-linear objective functions is more sophisticated as the PD update rule becomes non-linear in primal variables.
%Another body of research has studied the impact of quantized communications on the convergence behavior of consensus algorithms. In a consensus problem, a group of agents seek to arrive at a common state value. The authors in \cite{NOOT09} analyzed  the convergence time of average-consensus algorithms under quantization inter-agent communications and time-varying communication graphs. In \cite{KBS07}, the authors proposed a gossip-based algorithm for an average-consensus problem, and studied its convergence when agents' actions are integer-valued numbers. %Dong \emph{et al.} \cite{Dong13} analyzed the error performance of an average-consensus protocol, in a multi-agent consensus problem, wherein each agent, first,  quantizes its state value using a uniform quantizer, then, it transmits the quantized version of its state to other agents via a communication graph.
%Carli \emph{et al.} in \cite{CFSZ08} investigated the impacts of quantized inter-agent communications and communication topologies on the convergence rate of an average-consensus algorithm.  You and Xie \cite{YX11} derived necessary and sufficient conditions on the communication graph and the required  data rate for consensability of agents with linear local dynamics.

%We note that the impact of quantized communications on the convergence behavior of consensus algorithms has been extensively studied in the literature.  As this is not the main direction of our research, interested readers are referred to \cite{CFSZ08}-\cite{NOOT09} and references therein for a thorough investigation of this topic.
\subsection{Contributions}
%In this paper, we investigate the impact of quantized communications on the convergence behavior of the PD algorithm in a network utility maximization (NUM) problem. In our set-up,
We consider a NUM problem in which  $M$ agents  maximize the sum of their local concave objective functions subject to $N$ linear constraints using a quantized PD algorithm with a random initial condition. \textcolor{black}{As standard in the NUM literature, \emph{e.g.,} see  \cite{KMT98}, \cite{SS07} and references therein, we assume that the primal variables are updated by agents, and each dual variable is updated by a network node (NN) which has access to the knowledge of the the constraint associated with its dual variable}.% For example, routers in congestion control problems are embodiments of NNs.} %That is, each agent will update its corresponding primal variable using the knowledge of dual variables whereas the dual variables are updated by the NM using the knowledge of primal variables. %  a central entity called system which has only access to the constraint information. {\bf: Why constraints are known by system, why he does not know the objectives, give examples from resource allocation problems to justify} and each agent at each time updates its corresponding primal variable. The objective function of each agent depends on its local variable, however, the agents' {\bf decisions} are coupled through a set of constraints.  To update the primal and dual variables at each time, the system requires the knowledge of all primal variables in previous iteration, and the agents need the knowledge of all dual variables. This implies that, at each time,  agents have to transmit their decision variables to the system, and system needs to transmit the dual variables to agents at each time.
Thus, agents and NNs need to exchange the quantized values of primal and dual variables to execute the PD algorithm.% which are effectively the pricing signals. %However, in practice, the capacity of communication channels are limited by the distance between system and agents, available transmission power and bandwidth, etc.
%In practice, the communication between agents and NNs can only be performed in the quantized form due to the inherent capacity limitation of communication channels. This observation motivates us to study the convergence behavior of the PD algorithm in a more realistic scenario wherein agents and NNs can only exchange the quantized versions of primal and dual variables. 
We investigate the impact of quantized communications between agents and NNs on the rate of exponential mean-square convergence of the PD algorithm  under optimum achieving (OA) quantization schemes. The OA quantization schemes allow the primal and dual variables to converge to their optimal values as the time instance $k$ tends to infinity. 

First, using the information-theoretic notion of {\em differential entropy power}, we establish universal, explicit bounds on the fastest speed of asymptotic exponential mean square convergence for the PD, primal and dual variables to their corresponding optimal values (Theorem \ref{Theo: DDE}, \ref{Theo: DDE-P} and \ref{Theo: DDE-D}).
\textcolor{black}{Unlike previous studies of quantized optimization, a significant feature of these bounds is that they are completely independent of the OA quantization scheme employed,  making them 
applicable to {\em all} quantized PD algorithms. Given the utility functions, constraints and aggregate data rates (bits/sample)  of the agents and NN's, these results
give  system designers a way to determine in advance what exponential convergence speeds are  impossible to achieve.}
We note that the entropy power method has been used to study the stability of feedback control systems under quantization, \emph{e.g.,} see \cite{NE04} and \cite{FMS10}, \textcolor{black}{as well as convergence in quantized games \cite{NNA15}}.
%Here, our results show that the speed of exponential mean square convergence of the primal and dual variables depends on the average data rate used for transmission of dual and primal variables, respectively (see Theorems \ref{Theo: DDE-P} and \ref{Theo: DDE-D} for more details). % In the case of primal variables, the lower bound on the speed of exponential convergence is controlled by the number of primal variables, the behavior of second derivative of objective functions of agents around the optimal solution as well as the average aggregate data rate from the NM to agents (see Theorem \ref{Theo: DDE-P} for more details).  In the case of dual variables, our results show that the rate of exponential convergence of dual variables only depends on the number of constraints and the average aggregate data-rate from agents to the NM (see Theorem \ref{Theo: DDE-D}). This result suggests  that the convergence of dual variables is more sensitive to the communication data rate compared to that of primal variables. {\bf (Random initial condition)}

\textcolor{black}{Next, we obtain a bound on the fastest speed of exponential mean square convergence of PD variables in quadratic NUM problems under zoom-in quantization schemes (see Theorem \ref{Theo: EDE-New} for more details). This bound is significantly  tighter in the high data rate regime than Theorem \ref{Theo: DDE}.} We also derive  lower bounds on the mean square distance of PD, primal and dual variables from their corresponding optimal solutions for any given $k$ under quantized communication between agents and NNs (see Corollaries \ref{Coro: FTE} and \ref{Coro: FTE-P}). \textcolor{black}{Finally, we propose a uniform, zoom-in quantization scheme which allows the PD algorithm to converge to the optimal solution (Theorem 5).}

The organization of the paper is as follows. The next section describes our system model and assumptions. Section \ref{Sec: R&D} states our asymptotic and non-asymptotic results on the convergence of PD algorithm under quantization. In Section \ref{Sec: AOAQ}, we propose an OA quantization scheme. Section \ref{Sec: NR} presents our numerical results, and section \ref{Sec: Conc} concludes the paper. All  proofs are relegated to the appendices, to aid the fluency of the paper.

\section{System Model and Problem Formulation}\label{Sec: SM}
Consider a convex optimization problem in which $M$ agents maximize the sum of their individual objective functions subject to a set of linear equality constraints. % Let $\mathcal{M}=\left\{1,\cdots,M\right\}$ be the set of agents, and
 Let $x^i$ and $U_{i}\paren{x^i}$ represent the decision variable of agent $i$ and its objective function, respectively. It is assumed that the objective function of each agent is concave in its decision variable. The agents are interested in the solution of the following NUM problem:
  \begin{eqnarray}\label{Prob}
 \begin{array}{cc}
\underset{\vec{x}}{\rm maximize} & \sum\limits_{i}^M U_i\paren{x^i}\\
{\rm Subject\quad to} &  A\vec{x}= \vec{b}\\
  \end{array},
 \end{eqnarray}
where, $M$ is the number of agents, $\vec{b}\in \R^N$, $A\in \R^{N\times M}$, $N$ is the number of constraints,  and $\vec{x}=\left[x^1,\cdots,x^M\right]^\top$. We impose the condition $N<M$ to ensure that the feasible set of the optimization problem \eqref{Prob} is non-empty. The objective function in \eqref{Prob} is concave and the constraints are linear, thus, the optimization problem \eqref{Prob} can be solved using standard convex optimization techniques.

 Under the PD algorithm, the primal and dual variables are update according to
\begin{align}\label{EQ: NQU}
x^{i}_{k}&=x^{i}_{k-1}+\mu_{k-1} \paren{\frac{d }{d x^i}U_i\paren{x^i_{k-1}}-A^\top_i\vec{\lambda}_{k-1}},\nonumber  1\leq i\leq M\\
\lambda^{j}_{k}&=\lambda^j_{k-1}+\mu_{k-1}\paren{\bar{A}_j\vec{x}_{k-1}-b_j}\quad 1\leq j\leq N
\end{align}
, respectively, where $\mu_{k-1}$ is the step size of the algorithm at time $k-1$, $x^i_k$ and $\lambda^j_k$ denote the values of  $i$th primal variable and $j$th dual variable at time $k$, respectively, $\vec{\lambda}_{k-1}=\left[\lambda_{k-1}^1,\cdots,\lambda_{k-1}^N\right]^\top$, $A_i$ denotes the $i$th column of the matrix $A$ and $\bar{A}_j$ denotes the $j$th row of matrix $A$. To obtain the solution of the optimization problem \eqref{Prob}, we consider a primal-dual (PD) decomposition approach in which the primal variables, \emph{i.e.,} agents' decision variables, are updated by agents at each time. \textcolor{black}{Also, at each time step of the PD algorithm,  the $j$th dual variable, \emph{i.e.,} $\lambda^j$, is updated by the $j$th network node (NN) which has the knowledge of parameters characterizing the constraint associated with $\lambda^j$, \emph{i.e.,} $A_j$ and $b_j$.}  The vector of PD variables at time $k$, \emph{i.e.,} $\vec{y}_k$, is defined as the vector concatenation of $\vec{x}_k$ and $\vec{\lambda}_k$, \emph{i.e.,}
\begin{eqnarray}
\vec{y}_k=
\left[
\begin{array}{c}
\vec{x}_k\\
\vec{\lambda}_k
\end{array}
\right].\nonumber
\end{eqnarray}
%will employ the primal-dual algorithm. In each iteration of the primal-dual algorithm, the primal variables, \emph{i.e.,} decision variables, and the dual variables, \emph{i.e.,}  Lagrange multipliers, are simultaneously updated. Here, we assume that a separate entity called the system will update the dual variables {\bf: TODO: add justifications for the role of system}. Hence, at each time-step, agents will transmit the values of primal variables to the system, and system will inform agents about the values of the dual variables.
 In this paper, it is assumed that the initial primal and dual variables, \emph{i.e.,} $\vec{x}_0$ and $\vec{\lambda}_0$, are drawn randomly according to the probability density functions $p_{\vec{x}_0}\paren{\vec{x}}$ and $p_{\vec{\lambda}_0}\paren{\vec{\lambda}}$, respectively. \textcolor{black}{By allowing the initial condition to be random,  the primal and dual variables become random variables. This allows us to use information theoretic tools to study the speed of exponential convergence of the primal-dual algorithm under quantized communications.} We further impose the following assumptions on the objective functions of agents, step size $\mu_k$, $p_{\vec{x}_0}\paren{\vec{x}}$ and $p_{\vec{\lambda}_0}\paren{\vec{\lambda}}$. %{\bf: Complete the list of assumptions}
\begin{enumerate}
\item The agents' objective functions are concave and twice continuously differentiable.
\item $U^{\rm min}_i\leq \frac{d^2}{d {x^i}^2}U_i\paren{x^i}\leq U^{\rm max}_i<0$ for $x_i\in \R$ and all $i$.
\item $\mu_k\leq \min_i\frac{1}{\abs{U^{\rm min}_i}}$ for all $k$.
\item The sequence $\left\{\mu_k\right\}_k$ converges to $\mu^\star>0$.
\item The random vectors $\vec{x}_0$ and $\vec{\lambda}_0$ are mutually independent and the distributions of $\vec{x}_0$ and $\vec{\lambda}_0$ have finite differential entropies. That is,
\begin{eqnarray}
\left | - \int p_{\vec{x}_{0}}\paren{\vec{x}}\logp{p_{\vec{x}_{0}}\paren{\vec{x}}}d\vec{x}\right | <\infty\nonumber\\
\left | - \int p_{\vec{\lambda}_{0}}\paren{\vec{\lambda}}\logp{p_{\vec{\lambda}_{0}}\paren{\vec{\lambda}}}d\vec{\lambda}\right | <\infty\nonumber
\end{eqnarray}
\end{enumerate}
\textcolor{black}{Assumptions 1 and 2 are standard in the optimization literature. Assumption 2 implies that the objective functions of agents are strongly concave and the first derivative of each objective function is Lipschitz. Assumption 4 implies that the unquantized update rule does not employ a diminishing step-size rule as the PD update rule may not converge  exponentially with diminishing step-size rule. Assumptions 3 and 4, which are not commonly used in the literature, allow us to use the entropy power method. Assumption 5 implies that the initial condition injects a minimum  amount of uncertainty to the PD algorithm, and the amount of uncertainty due to the initial condition is bounded. Variants of assumption 5 are used in the quantized feedback control literature, \emph{e.g.,} see \cite{NE04}.
%EHSAN:  PLEASE ALSO INDICATE WHICH OF THESE ASSUMPTIONS IS STANDARD/NONSTANDARD IN OPTIMIZATION, AS REQUESTED BY REV. 3 
}

\subsection{Quantizer Structure}

To execute the PD update rule \eqref{EQ: NQU}, the agents and NNs require the knowledge of dual and primal variables, respectively. 
Since the agents and NNs are not necessarily co-located, the information exchange between NNs and agents is performed via broadcast communication channels, as described in the next subsection. 
Due to the capacity limitations of these channels,  only quantized versions of the primal and dual variables can be exchanged between NNs and agents. %deployed between agents and the NM, agents and the NM can only exchange the quantized versions of their variables which are representable using finite number of bits.

\textcolor{black}{To transmit $x^i_k$ to NNs, agent $i$ encodes $x^i_k$ to $\hat{Q}^{\vec{x}}_{i,k}$ using an adaptive encoder mapping  of the form
\begin{eqnarray}
\hat{Q}^{\vec{x}}_{i,k}=E_{i,k}^{\vec{x}}\paren{\left\{x^i_n\right\}_{n=0}^k,\left\{\hat{Q}^{\vec{x}}_{i,n}\right\}_{n=0}^{k-1}}. \nonumber
\end{eqnarray} 
 Then, it broadcasts $\hat{Q}^{\vec{x}}_{i,k}$ to NNs. The output of the encoder of agent $i$ at time $k$, \emph{i.e.,} $\hat{Q}^{\vec{x}}_{i,k}$, belongs to the finite alphabet set $\mathcal{A}^{\vec{x}}_{i,k}$. Thus, agent $i$ requires $\log_2\abs{\mathcal{A}^{\vec{x}}_{i,k}}$ bits to transmit its encoded symbol to NNs.    A large value of $\abs{\mathcal{A}^{\vec{x}}_{i,k}}$ indicates that agent $i$ transmits its decision variable with a high precision to NNs whereas a low $\abs{\mathcal{A}^{\vec{x}}_{i,t}}$  indicates low quality communication between agent $i$ and NNs.}

 \textcolor{black}{Upon receiving $\hat{Q}^{\vec{x}}_{i,k}$, each NN reconstructs a  quantized estimate of $x^i_k$, \emph{i.e.,} $\quant{}{\vec{x}}{i,k}$, 
using the decoder mapping $\quant{}{\vec{x}}{i,k}=D^{\vec{x}}_{i,k}\paren{\left\{\hat{Q}^{\vec{x}}_{i,n}\right\}_{n=0}^k}$.
}
\textcolor{black}{Similarly, at time $k$, NN $j$ chooses the symbol $\hat{Q}^{\vec{\lambda}}_{j,k}$ from the finite alphabet set $\mathcal{A}^{\vec{\lambda}}_{j,k}$ according to the adaptive encoding map 
\begin{eqnarray}\hat{Q}^{\vec{\lambda}}_{j,k}=E_k^{\lambda}\paren{\left\{\lambda^j_n\right\}_{n=0}^k,\left\{\hat{Q}^{\vec{\lambda}}_{j,n}\right\}_{n=0}^{k-1}},\nonumber
\end{eqnarray}
 and broadcasts $\hat{Q}^{\vec{\lambda}}_{j,k}$ to all agents. Then, agents construct the quantized version of $\lambda^j_{k}$, \emph{i.e.,} $\quant{}{\vec{\lambda}}{j,k}$, using the decoding map $\quant{}{\vec{\lambda}}{j,k}=D^{\vec{\lambda}}_{j,k}\paren{\left\{\hat{Q}^{\vec{\lambda}}_{j,n}\right\}_{n=0}^k}$. Note that our formulation allows the encoded symbol at time $k$ to depend on the current and past values of primal/dual variables as well as the past outputs of the encoder.}
 \textcolor{black}{We refer to $\mathcal{Q}=\!\!\left\{\!\!\left\{\!E^{\vec{x}}_{i,k}\paren{\cdot}\!,\!D^{\vec{x}}_{i,k}\paren{\cdot}\!\right\}_i\!,\!\left\{\!E^{\vec{\lambda}}_{j,k}\paren{\cdot}\!, \!D^{\vec{\lambda}}_{j,k}\paren{\cdot}\!\right\}_j\!\!\right\}_{k=0}^\infty$ as a quantization scheme.}
Also, the quantized version of the PD variables at time $k$ under the quantization scheme $\mathcal{Q}$ is denoted by $\quant{}{}{k}$, \emph{i.e.,}
\begin{eqnarray}
\quant{}{}{k}=
\left[
\begin{array}{c}
\quant{}{\vec{x}}{k}\\
\quant{}{\vec{\lambda}}{k}\nonumber
\end{array}
\right],
\end{eqnarray}
where $\quant{}{\vec{x}}{k}\!\!=\!\!\left[\quant{}{\vec{x}}{1,k},\cdots,\quant{}{\vec{x}}{M,k}\right]^\top$ and $\quant{}{\vec{\lambda}}{k}\!\!=\!\!\left[\quant{}{\vec{\lambda}}{1,k},\cdots,\quant{}{\vec{\lambda}}{N,k}\right]^\top$. 

Next, we define three notions of data rate for a given quantization scheme $\mathcal{Q}$. Later, these data rates are used to study the convergence behavior of primal, dual and PD variables. The average aggregate data rate per unit time for transmitting the primal variables to NNs under the quantization scheme $\mathcal{Q}$, $R_{\vec{x}}$, is defined as
\begin{eqnarray}
R_{\vec{x}}=\lim_{k\rightarrow\infty}\frac{1}{k}\sum_{t=0}^{k-1}\paren{\sum_{i=1}^M\log\abs{\mathcal{A}^{\vec{x}}_{i,t}}}
\end{eqnarray}
Similarly, we define the average aggregate data rate per unit time for broadcasting the dual variables to agents under the quantization scheme $\mathcal{Q}$, $R_{\vec{\lambda}}$, as
\begin{eqnarray}\label{Eq: R-dual}
R_{\vec{\lambda}}=\lim_{k\rightarrow\infty}\frac{1}{k}\sum_{t=0}^{k-1}\paren{\sum_{j=1}^N\log\abs{\mathcal{A}^{\vec{\lambda}}_{j,t}}}
\end{eqnarray}
 Finally, the average total date rate per unit time under the quantization scheme $\mathcal{Q}$, \emph{i.e.,} $R_{\mathcal{Q}}$, is defined as
\begin{eqnarray}\label{Eq: R_total}
R_{\mathcal{Q}}=\lim_{k\rightarrow\infty}\frac{1}{k}\sum_{t=0}^{k-1}\paren{\paren{\sum_{i=1}^M\log\abs{\mathcal{A}^{\vec{x}}_{i,t}}}+\sum_{j=1}^N\log\abs{\mathcal{A}^{\vec{\lambda}}_{j,t}}}
\end{eqnarray}

The quantized PD update rule under the quantization scheme $\mathcal{Q}$ can be written as
\begin{align}\label{EQ: QU}
x^{i}_{k}&=x^{i}_{k-1}+\mu_{k-1} \paren{\frac{d }{d x^i}U_i\paren{x^i_{k-1}}-A^\top_i\quant{\vec{\lambda}_{k-1}}{\vec{\lambda}}{k-1}},\nonumber\\
\lambda^{j}_{k}&=\lambda^j_{k-1}+\mu_{k-1}\paren{\bar{A}_j\quant{\vec{x}_{k-1}}{\vec{x}}{k-1}-b_j}
\end{align}
%Using the quantizer introduces a distortion to the primal and dual variable. That is, the primal and dual variables of the optimization algorithm under quantization are different from the primal and dual variables in the non-quantized algorithm.
Let $\vec{x}^\star$, $\vec{\lambda}^\star$ be the optimal primal and dual solutions, respectively. Also, let  $\vec{y}^\star$ be the vector concatenation of  $\vec{x}^\star$, $\vec{\lambda}^\star$. We define $\epsilon_k=\vec{y}_k-\vec{y}^\star$ as the difference between the PD variables at time $k$ and the optimal solution. Let $\norm{\vec{\epsilon}_k}{2}$ denote the distance of the PD variables at time $k$ from optimal solution, \emph{i.e.,}
\begin{eqnarray}
\norm{\vec{\epsilon}_k}{2}=\sqrt{\sum_{i=1}^M\paren{x^i_k-{x^i}^\star}^2+\sum_{j=1}^N\paren{\lambda^j_k-{\lambda^j}^\star}^2}
\end{eqnarray}
where ${x^i}^\star$ and ${\lambda^j}^\star$ are the optimal values of the primal variable $x^i$ and the dual variable $\lambda^j$, respectively. Then, the mean square distance (MSD) of the PD variables from the optimal solution at time $k$ under the quantization scheme $\mathcal{Q}$ is defined as $\ES{\norm{\vec{\epsilon_k}}{2}^2}$. We define the MSD of the primal variables from the optimal primal solution at time $k$ as $\ES{\norm{\vec{\epsilon}_k^{\vec{x}}}{2}^2}$ where $\vec{\epsilon}_k^{\vec{x}}=\vec{x}_k-\vec{x}^\star$. Similarly, the MSD of dual variables at time $k$ from the optimal dual solution is defined as $\ES{\norm{\vec{\epsilon}_k^{\vec{\lambda}}}{2}^2}$ where $\vec{\epsilon}_k^{\vec{\lambda}}=\vec{\lambda}_k-\vec{\lambda}^\star$. %In the next section,  be the difference between the primal variables and their optimal solution, respectively, defined as
 %In this paper, we establish fundamental connections between the decay rate of the distortion and the the average rate per iteration for different decomposition architectures of the optimization problem \eqref{Prob}.
Next, we define the class of optimum achieving (OA) quantization schemes.
\begin{definition}\label{Def: OA}
The quantization scheme $\mathcal{Q}$ is called an OA quantization scheme if, under $\mathcal{Q}$, the primal and dual variables converge to their optimal values $\vec{x}^\star$ and $\vec{\lambda}^\star$, receptively. That is: %{\bf: is the second condition necessary?}
\begin{eqnarray}
\lim_{k\rightarrow\infty}\vec{x}_k=\vec{x}^\star\nonumber\\
\lim_{k\rightarrow\infty}\vec{\lambda}_k=\vec{\lambda}^\star\nonumber
\end{eqnarray}
\end{definition}
The Definition \ref{Def: OA} implies that, under an OA quantization scheme, the quantization error does not impede the convergence of the  PD algorithm to the optimal solution. Thus, under an OA quantization scheme, the PD algorithm converges to the optimal solution of the optimization problem regardless of the quantized  communication between agents and NNs.

\subsection{Communication Graph and Communication Cost}

%we extend our results to a more general network, induced by the $N\times M$ constraint matrix $A$.

%Note that in typical distributed problems,  $A$ may only have a relatively small number of nonzero entries,
%reflecting for instance neighbouring agents that share some common resource.

\textcolor{black}{The communication topology is represented  by a bipartite graph induced by the $N\times M$ constraint matrix $A$.  In this graph, edges exist only between agents and  network nodes (NNs), which  form two disjoint sets of vertices.
There exists an edge  between agent $i$ and NN $j$ 
  in the communication graph if and only if $A_{ji}\neq 0$. The communication mechanism is broadcast in nature,
with each vertex  `listening' and broadcasting only to those other vertices with which it shares an edge.
This is implemented by uniquely assigning every vertex in the graph   one of $N+M$ disjoint transmission radio-frequency bands ({\em frequency division multiplexing}) or one of $N+M$ disjoint time slots per cycle ({\em time division multiplexing}),
before the system is deployed. Any other vertex that needs to listen to a transmission just tunes in to the appropriate frequency band or time slot dedicated to the corresponding transmitter. 
Note that the edges do not represent individual one-to-one channels, but indicate the  broadcast transmitter-receiver structure of the system.}

\textcolor{black}{Under typical digital modulation formats, the width of the frequency band/time-slot allocated to agent $i$ and/or the average transmission power it consumes to broadcast its encoded  symbols 
to all NNs $j$ with $A_{ji}\neq 0$ will be proportional to its average data rate $R^i_{\vec{x}}:=\lim_{k\rightarrow\infty}\frac{1}{k}\sum_{t=0}^{k-1}\log\abs{\mathcal{A}^{\vec{x}}_{i,t}}$. 
Similarly, the band/slot-width and/or transmission power used by NN $j$ to broadcast its encoded dual symbols 
to all agents $i$ with $A_{ji}\neq0$  is typically proportional to $R^j_{\vec{\lambda}}:=\lim_{k\rightarrow\infty}\frac{1}{k}\sum_{t=0}^{k-1}\log\abs{\mathcal{A}^{\vec{\lambda}}_{j,t}}$.
Equation \eqref{Eq: R_total}, which can be intuitively interpreted as $\sum_{i=1}^M R^i_{\vec{x}}+\sum_{j=1}^NR^j_{\vec{\lambda}}$, then captures the total amount of physical resources, \emph{i.e.,} time, bandwidth or power,  required for the system to communicate. It can be seen that this cost scales like $O(N+M)$ as the network grows in size. Note that due to the broadcast nature of the system,
 every transmission can be heard by multiple receivers, without the transmitter having to use up extra resources. 
%EHSAN - THIS JUSTIFICATION WORKS ONLY IF THE LIMSUP = THE STRAIGHT LIM.
%HOWEVER,  IT IS A VERY MILD ASSUMPTION THAT THE STRAIGHT LIMIT EXISTS IN THE DATA RATE DEFN.
%SO WHY DON;T WE CHANGE THE DATA RATE DEFN TO ASSUME THAT THE STRAIGHT LIMIT EXISTS
}

%\textcolor{black}{The total number of transmitted bits at time $k$ to NNs by agents is $\sum_{i=1}^M\log\abs{\mathcal{A}^{\vec{x}}_{i,k}}$. Hence, the average number of bits per agent $\frac{1}{M}\sum_{i=1}^M\log\abs{\mathcal{A}^{\vec{x}}_{i,k}}$ remains %bounded as $M$ becomes large. Similarly, the total number of bits used by NNs to transmit the encoded dual variables at time $k$ is equal to $\sum_{j=1}^N \log{\abs{\mathcal{A}^{\vec{\lambda}}_{j,k}}}$ which is independent of number of agents. We note that the %communication channel between each agent and NNs is a broadcast channel with \emph{only} common message as each agent transmits a common message to NNs. Similarly, the communication channel between each NN and agents is in the form of a broadcast channel %with only common message.}

\section{Results and Discussions}\label{Sec: R&D}
In this section, we analyze the impact of quantized communications on the mean square distance (MSD) of primal-dual (PD), primal and dual variables from the optimal solution in two different regimes: $(i)$ Asymptotic regime, $(ii)$ Non-asymptotic regime. In the asymptotic regime, we are concerned with the behavior of MSD under OA quantization schemes as the time index $k$ tends to infinity. To this end, the notion of distance decay exponent (DDE) is introduced which captures the rate of exponential convergence of MSD to zero. % mean square convergence of the PD algorithm to the optimal solution under an OA quantization scheme.
We establish universal lower bounds on the DDE of PD variables, primal variables and dual variables (see Theorems \ref{Theo: DDE}, \ref{Theo: DDE-P}, \ref{Theo: DDE-D} and \ref{Theo: EDE-New} for more details).
In the non-asymptotic regime, we are concerned with the behavior of the MSD for any finite $k$. Here, our results provide universal lower bounds on the MSD of PD, primal and dual variables, from the optimal solution, for any finite $k$ (see Corollaries \ref{Coro: FTE} and \ref{Coro: FTE-P} for more details). We start by presenting our asymptotic results in the next subsection.
\subsection{Asymptotic behavior of MSD in PD algorithm}
In this subsection, first, we introduce the notion of distance decay exponent (DDE) for the PD, primal and dual variables.  Then, we derive universal lower bounds on the DDE of PD, primal and dual variables.
\begin{definition}
Let $\mathcal{Q}$ be an OA quantization scheme. Then, the DDE of the PD, primal and dual variables under $\mathcal{Q}$ are  defined as
\begin{eqnarray}
\liminf_{k\rightarrow\infty}\frac{1}{k}\log\ES{\norm{\vec{\epsilon}_k}{2}^2},\nonumber\\
\liminf_{k\rightarrow\infty}\frac{1}{k}\log\ES{\norm{\vec{\epsilon}_k^{\vec{x}}}{2}^2},\nonumber\\
\liminf_{k\rightarrow\infty}\frac{1}{k}\log\ES{\norm{\vec{\epsilon}_k^{\lambda}}{2}^2},\nonumber
\end{eqnarray}
respectively.
\end{definition}
The DDEs capture the  speed of exponential mean square convergence of the PD, primal and dual variables to their corresponding optimal solutions. They are non-positive quantities, and a more negative DDE indicates faster convergence to the optimal solution. \textcolor{black}{Also, a zero DDE implies slower-than-exponential convergence.} In this subsection, the information-theoretic notion of entropy power is used to establish universal lower bounds on the DDE of the PD/primal/dual variables. % The lower bound is universal in the sense that it is independent of the structure of quantizer.

The next theorem provides a universal lower bound on the  DDE of the PD variables under OA quantization schemes.
\begin{theorem}\label{Theo: DDE}
Let  $\mathcal{Q}$  be an OA quantization scheme. Then, the DDE of PD variables under $\mathcal{Q}$ can be lower bounded as
\begin{align}
&\liminf_{k\rightarrow\infty}\frac{1}{k}\log\ES{\norm{\vec{\epsilon}_k}{2}^2}\nonumber\\
&\hspace{1cm}\geq\frac{2}{N+M}\paren{\sum_{i=1}^M\log\paren{1+\mu^\star \frac{d^2}{d {x^i}^2}U_i\paren{{x^i}^\star}}-R_{\mathcal{Q}}}.
\end{align}
where ${x^i}^\star$ is the optimal value of the primal variable $x^i$.
\end{theorem}
\begin{IEEEproof}
Please see Appendix \ref{App: DDE}.
\end{IEEEproof}
Theorem \ref{Theo: DDE} establishes an explicit universal lower bound on the DDE of PD variables under OA quantization schemes. 
This bound is universal in the sense that it is independent of the structure of quantizer, and is thus applicable to all quantization schemes.

According to Theorem \ref{Theo: DDE}, for a given average total data rate $R_{\mathcal{Q}}$, the PD variables  converge to the optimal solution at most exponentially fast.
The speed of this exponential convergence is bounded by the average total data rate under the quantization scheme, \emph{i.e.,} $R_{\mathcal{Q}}$, 
and also by the behavior of the objective functions of agents around the optimal solution. As stated in Theorem \ref{Theo: DDE}, the lower bound on the DDE for PD variables decreases linearly with $R_{\mathcal{Q}}$.
Note that as $R_{\mathcal{Q}}$ becomes large, the NNs and agents have more precise information about the primal and dual variables. \textcolor{black}{Thus one might intuitively expect a  quantized PD algorithm to converge faster 
to the optimal solution as $R_{\mathcal{Q}}$ increases. The result above is consistent with this intuition.}

The lower bound on the DDE also increases with the second derivatives of the agents' objective functions  at the optimal solution. 
As these second derivatives becomes less negative, 
the objective function becomes flatter near the optimal solution \textcolor{black}{and the quantized PD algorithm can be expected to converge more slowly. This result is also in concordance with this intuition.}

The next theorem establishes a universal lower bound on the DDE of primal variables in the quantized PD update rule under an OA quantization scheme.
\begin{theorem}\label{Theo: DDE-P}
Consider the OA quantization scheme $\mathcal{Q}$. Then, the DDE of the primal variables under $\mathcal{Q}$ is lower bounded as
\begin{align}
&\liminf_{k\rightarrow\infty}\frac{1}{k}\log\ES{\norm{\vec{\epsilon}^{\vec{x}}_k}{2}^2}\nonumber\\
&\hspace{1cm}\geq\frac{2}{M}\paren{\sum_{i=1}^M\log\paren{1+\mu^\star \frac{d^2}{d {x^i}^2}U_i\paren{{x^i}^\star}}-R_{\vec{\lambda}}}.
\end{align}
\end{theorem}
\begin{IEEEproof}
Please see Appendix \ref{App: DDE-P}.
\end{IEEEproof}
According to Theorem \ref{Theo: DDE-P},  the  exponential convergence speed of the primal variables is limited by the behavior of objective functions of agents around the optimal solution, the average aggregate data rate for transmission of dual variables and the number of agents. Different from the PD bound in Theorem \ref{Theo: DDE}), this  lower bound on the DDE of the primal variables depends only on the average aggregate data rate for transmission of dual variables, \emph{i.e.,} $R_{\vec{\lambda}}$, rather than on the average total data rate under the quantization scheme $\mathcal{Q}$. This observation signifies the role of the quantized dual variables on the convergence of the primal variables. %Note that as $R_{\vec{\lambda}}$ increases, the quantization error of dual variables reduces. Thus, agents update the primal variables using more accurate values of dual variables which suggests a faster convergence to the optimal solution.

In the next theorem, we study the DDE for dual variables.
\begin{theorem}\label{Theo: DDE-D}
The DDE of dual variables under the OA quantization scheme $\mathcal{Q}$ satisfies 
\begin{align}
&\liminf_{k\rightarrow\infty}\frac{1}{k}\log\ES{\norm{\vec{\epsilon}^{\vec{\lambda}}_k}{2}^2}\geq-\frac{2}{N}R_{\vec{x}}.
\end{align}
\end{theorem}
\begin{IEEEproof}
Please see Appendix \ref{App: DDE-D}.
\end{IEEEproof}
Theorem \ref{Theo: DDE-D} establishes a universal  bound on the fastest possible exponential convergence rate of the dual variables under any OA quantization  scheme $\mathcal{Q}$. The lower bound in Theorem \ref{Theo: DDE-D} is controlled by the number of constraints and the average aggregate data rate for transmission of primal variables to NNs. Compared to the PD lower bound, it does not depend on the behavior of the objective functions of agents and is only limited by the average aggregate data rate for transmission of the primal variables, \emph{i.e.,} $R_{\vec{x}}$, rather than the average total data rate  $R_{\mathcal{Q}}$. %Since the update rule for dual variables is linear with respect to the dual variables, the lower bound in Theorem \ref{Theo: DDE-D} is independent of objective functions of agents.
%According to Theorem \ref{Theo: DDE-D}, for a fixed $R_{\vec{x}}$, the lower bound on the DDE of dual variables increases as the number of constraints becomes large. Thus, to maintain the same lower bound on the DDE of dual variables, the average aggregate data rate for transmission of primal variables should increase as  $N$ increases.

\textcolor{black}{Next, we derive a lower bound on the DDE of the PD algorithm in quadratic NUM problems under zoom-in quantization schemes (see Definition \ref{Def: Zoom-in}). This bound is tighter compared with the lower bound in Theorem \ref{Theo: DDE} at the high data-rate regime. In a quadratic NUM problem, the objective function of agent $i$ is given by $U_i\paren{x^i}=-\frac{a_i}{2}\paren{x^i}^2+c_ix^i+f_i$ where $a_i$ is a positive constant. The unquantized PD algorithm for quadratic NUM problems can be written as  
\begin{align}\label{Eq: QUR}
x^{i}_{k}&=\paren{1-\mu a_i}x^{i}_{k-1}+\mu \paren{c_i-A^\top_i\vec{\lambda}_{k-1}},\nonumber  1\leq i\leq M\\
\lambda^{j}_{k}&=\lambda^j_{k-1}+\mu\paren{\bar{A}_j\vec{x}_{k-1}-b_j}\quad 1\leq j\leq N
\end{align}
 Let $\vec{y}_k$ be the vector concatenation of $\vec{x}_k$ and $\vec{\lambda}_k$. Then, \eqref{Eq: QUR} can be written as 
\begin{align}
\vec{y}_{k}=T\vec{y}_{k-1}+\mu
\left[
\begin{array}{c}
\vec{c}\\
-\vec{b}
\end{array}
\right]\nonumber
\end{align}
 where $\vec{c}=\left[c_1\cdots,c_M\right]^\top$ and the matrix $T$ is defined as  
\begin{align}
T=\left[
 \begin{array}{cc}
	{\rm Diag}\paren{1-\mu a_1,\cdots,1-\mu a_M}&-\mu A^\top\\
	\mu A&I_{N}
\end{array}
\right]
\end{align}
in which $I_N$ denotes an $N$-by-$N$ identity matrix and ${\rm Diag}\paren{1-\mu a_1,\cdots,1-\mu a_M}$ is a diagonal matrix with the $i$th diagonal element equal to $1-\mu a_i$.} 

\textcolor{black}{Let $\tilde{Q}_k=\left\{\hat{Q}^{\vec{x}}_{1,n},\cdots,\hat{Q}^{\vec{x}}_{M,n},\hat{Q}^{\vec{\lambda}}_{1,n},\cdots,\hat{Q}^{\vec{\lambda}}_{N,n}\right\}_{n=0}^k$ be the collection of encoders' outputs up to time $k$, respectively. The quantized PD update rule is denoted by $\vec{y}_{k+1}=\QM{\vec{y}_k}{\tilde{q}_k}$ where $\tilde{q}_k$ is a realization of $\tilde{Q}_k$. We use $C_k\paren{\tilde{q}_k}$ to represent the quantization cell corresponding to $\tilde{q}_k$. Next, a zoom-in quantization scheme is defined.}
\begin{definition}\label{Def: Zoom-in}
 \textcolor{black}{Consider the quantization scheme $\mathcal{Q}$, and let $C_k\paren{\tilde{q}_k}$ be the quantization cell at time $k$  which contains $\vec{y}_k$. Then, $\mathcal{Q}$ is a zoom-in quantization scheme if at time $k+1$ the image of $C_k\paren{\tilde{q}_k}$ under $\QM{\cdot}{\tilde{q}_k}$ is quantized for all $k\in\N_0=\left\{0,1,2,\cdots\right\}$.}
\end{definition}
%\begin{figure}[!t]
%\centering{\includegraphics[scale=0.4]{Image.eps}}
%\caption{Two dimensional lattice of integers $\Z^2$ $(a)$ and the lattice $T\Z^2$ (b).} \label{F-Image}
%\end{figure}

\textcolor{black}{In addition to the assumptions in Section \ref{Sec: SM}, we assume that 
\begin{enumerate}
\item The matrix $T$ is invertible and all its eigenvalues are inside the unit circle in complex plane. 
\item A zoom-in quantization scheme is employed and each primal/dual variable is independently quantized. 
\item The distributions of initial primal and dual variables, \emph{i.e.,} $p_{\vec{x}_0}\paren{\vec{x}}$ and $p_{\vec{\lambda}_0}\paren{\vec{\lambda}}$, are bounded and have  finite support sets.
\end{enumerate} }

\begin{figure}[!t]
\centering{\includegraphics[scale=0.5]{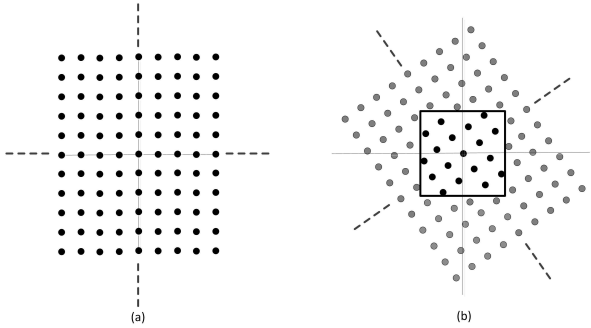}}
\caption{Two dimensional lattice of integers $\Z^2$ $(a)$ and the lattice $T\Z^2$ (b).} \label{F-lattice}
\end{figure}

\begin{theorem}\label{Theo: EDE-New}
\textcolor{black}{Consider  any zoom-in quantization scheme $\mathcal{Q}$ with $\rho=\frac{\delta^{\rm max}_k}{\delta^{\rm min}_k}$ (for all $k$) where $\delta^{\rm max}_k$ and $\delta^{\rm min}_k$ are the maximum and minimum quantization steps under $\mathcal{Q}$ at time $k$, respectively. Let ${\rm B}$ be the hypercube centered at the origin with the $i$th side length equal to $4\rho\abs{T_{ii}}+2\norm{T}{\infty}$ where $\norm{\cdot}{\infty}$ denotes the norm infinity and $T_{ii}$ is the $i$ diagonal entry of matrix $T$. Let $\beta_T$ be the number elements in the set ${\rm B}\cap T\Z^{N+M}$ where the lattice $T\Z^{N+M}$ is defined as $T\Z^{N+M}=\left\{T\vec{I},\quad \vec{I}\in\Z^{N+M}\right\}$ and $\Z^{N+M}$ is the lattice of integers in $\R^{N+M}$. %Also, let  $\beta_T$ be the number of integer lattice points $\vec{I}\in\Z^{M+N}$ which satisfy the following system of inequalities 
%\begin{align}
%\abs{T\vec{I}}\leq 2\rho T_d \nonumber
%\end{align}
%where $T_d=\left[\abs{T_{ii}}+\frac{\norm{T}{\infty}}{2\rho}\right]_i^\top$, $\norm{\cdot}{\infty}$ denotes the norm infinity, $T_{ii}$ is the $(i,i)$ entry of matrix $T$ and $\leq$ indicates component-wise inequality. 
Then, the DDE of the PD variables  under $\mathcal{Q}$ for quadratic NUM problems is lower bounded as 
\begin{align}\label{Eq: EDE-New}
\liminf_{k\rightarrow\infty}\frac{1}{k+1}&\log{\ES{\norm{\vec{\epsilon}_{k+1}}{2}^2}}\nonumber\\
&\geq -\frac{2}{M+N}\logp{\frac{ \beta_T}{\paren{\prod_{i=1}^{M+N}\abs{T_{ii}}}}} 
\end{align}.}
\end{theorem}
\begin{IEEEproof}
\textcolor{black}{Please see Appendix \ref{App: EDE-New}.}
\end{IEEEproof}
\textcolor{black}{Theorem \ref{Theo: EDE-New} establishes a bound on the fastest possible exponential convergence speed of quantized PD algorithms in quadratic NUM problems, under any zoom-in quantization scheme. The lower bound in Theorem \ref{Theo: EDE-New} depends on the number of agents, number of constrains and $\beta_T$. The constant $\beta_T$ depends on the dynamics of unquantized PD algorithm, \emph{i.e.,} matrix $T$, and can be interpreted as the number of lattice points in $\Z^{N+M}$ which lie in ${\rm B}$ after applying the linear transformation $T$ to $\Z^{N+M}$. Fig. \ref{F-lattice} shows the two dimensional lattice of integers $\Z^2$ and its image after applying a linear transformation. In Fig. \ref{F-lattice} $(b)$, the number of lattice points in the square is equal to $\beta_T$. Since the transformation $T$ is linear, $\vec{0}$ always lies in ${\rm B}$ which implies $\beta_T\geq 1$.}

\textcolor{black}{Consider the PD algorithm in a quadratic NUM problem under the zoom-in quantization scheme $\mathcal{Q}$ with $\rho=\frac{\delta^{\max}_k}{\delta^{\min}_k}$. For the PD algorithms, Theorems \ref{Theo: DDE} and \ref{Theo: EDE-New} can be combined into
\begin{align}\label{Eq: LB-COMB}
&\liminf_{k\rightarrow\infty}\frac{1}{k+1}\log{\ES{\norm{\vec{\epsilon}_{k+1}}{2}^2}}\geq \nonumber\\
&\frac{2}{M+N}\paren{\paren{\sum_{i=1}^M \logp{1-\mu a_i}}-\min\paren{\logp{\beta_T},R_{\mathcal{Q}}} }
\end{align}
 If the quantization intervals for each primal/dual variable is divided to $K\geq 2$ equal length intervals, data-rate under the quantization $\mathcal{Q}$ \emph{i.e.,} $R_{\mathcal{Q}}$,  will increase by $\paren{N+M}\logp{K}$ bits and $\rho$ does not change. Hence, according to \eqref{Eq: LB-COMB}, the lower bound in Theorem \ref{Theo: EDE-New} becomes tighter compared to that in Theorem \ref{Theo: DDE} as $R_{\mathcal{Q}}$ (or $K$) becomes large. This observation shows that the exponential convergence speed of the quantized PD algorithm in quadratic NUM problems cannot be made arbitrarily fast by increasing $R_{\mathcal{Q}}$. }

\textcolor{black}{An upper bound on $\beta_T$ can be obtained by finding the number of lattice points of $\Z^{N+M}$ which lie in the smallest hypercube containing the image of ${\rm B}$ under $T^{-1}$. Let $T^{-1}\paren{\rm B}$ be the image of the hypercube ${\rm B}$ under linear transformation $T^{-1}$. Let ${\rm B}^\star_{T^{-1}}$ be the smallest hypercube containing $T^{-1}\paren{\rm B}$. Then, $\beta_T$ is upper bounded by $\prod_{i}\paren{\left\lfloor l^\star_i\right\rfloor+1}$ where $l^\star_i$ is the $i$th side length of ${\rm B}^\star_{T^{-1}}$. In our numerical results, this upper bound on $\beta_T$ is used to compute the lower bound in Theorem \ref{Theo: EDE-New}.}
\subsection{MSD of the PD algorithm in non-asymptotic regime}
In this subsection, we establish universal lower bounds on the mean square distance (MSD) of primal-dual (PD), primal and dual variables from their corresponding optimal solutions at any finite time instance $k$. Unlike Theorems \ref{Theo: DDE}, \ref{Theo: DDE-P} and \ref{Theo: DDE-D}, the following results are not limited to optimum achieving (OA) quantization schemes. Thus, they give rise to universal lower bounds on the MSD of PD, primal and dual variables from their corresponding optimal solutions, under arbitrary quantization schemes. Our results in this subsection indicate that the distance between the optimization variables and the  optimal solution cannot be made arbitrarily close to zero at a given time instance $k$.  %They also provide lower bounds on the closeness of the optimization variables to their optimal solutions at any given time instance under the quantized PD algorithm. 
We start by presenting our non-asymptotic lower bound on the MSD of the PD variables.
\begin{corollary}\label{Coro: FTE}
Consider the PD algorithm under the quantization scheme $\mathcal{Q}$. Then, the MSD of the PD variables from the optimal solution at time $k$ can be lower bounded as
\begin{align}\label{Eq: LB-G}
&\log\ES{\norm{\vec{\epsilon}_k}{2}^2}\geq \logp{\frac{\e{1-\frac{1}{M+N}}}{{2\pi \e{}}}}+ \nonumber\\
&\hspace{1cm}\frac{2}{N+M}\left(\sum_{i=1}^M\sum_{n=0}^{k-1}{\logp{1+\mu_nU^{\min}_i}}\right.\nonumber\\
&\hspace{1cm}+\left.\mathsf{h}\left[\vec{y}_{0}\right]-\sum_{t=0}^{k-1}\paren{\paren{\sum_{i=1}^M\log\abs{\mathcal{A}^{\vec{x}}_{i,t}}}+\sum_{j=1}^N\log\abs{\mathcal{A}^{\vec{\lambda}}_{j,t}}}\right),
\end{align}
\end{corollary}
\begin{IEEEproof}
The proof directly follows from inequalities \eqref{Eq: MS-LB} and \eqref{Eq: Aux-1} in Appendix \ref{App: DDE}.
\end{IEEEproof}
Corollary \ref{Coro: FTE} provides a universal lower bound on the MSD of PD variables under quantized communications between agents and NNs.  This result indicates that at a given time the PD variables cannot be arbitrarily close to the optimal solution (in the mean square sense), and imposes a lower bound on the MSD of PD variables from the optimal solution at a given time. According to Corollary \ref{Coro: FTE}, the MSD of PD variables from the optimal solution at time $k$ is bounded from below by the behavior of second derivative of objective functions of agents along the trajectories of primal variables up to time $k-1$, the total number of bits exchanged between agents and NNs up to time $k-1$, the differential entropy of distribution of initial PD variables, \emph{i.e.,} $\mathsf{h}\left[\vec{y}_{0}\right]$, and the number of constraints and agents. The impacts of objective functions of agents and the data rate between agents and NNs on the lower bound in \eqref{Eq: LB-G} are similar to those in Theorem \ref{Theo: DDE}, and they are not discussed to avoid repetition.

Note that the entropy power of $\vec{y}_{0}$, \emph{i.e.,} $\frac{1}{2\pi\e{}}\e{\frac{2}{N+M}\mathsf{h}\left[\vec{y}_{0}\right]}$ is a measure of effective support volume of the random vector $\vec{y}_{0}$. Thus, as $\mathsf{h}\left[\vec{y}_{0}\right]$ becomes large, the size of effective support set of $\vec{y}_{0}$ increases, \emph{i.e.,} $\vec{y}_0$ will be distributed on a larger region of $\R^{N+M}$. As a result, the MSD of the PD variables from the optimal solution increases since $\vec{y}_{0}$ effectively takes value from a larger set,  a behavior predicted by Corollary \ref{Coro: FTE}.

%   on the MSD, assume that $\vec{y}_{0}$ is Gaussian distributed with zero mean and covariance matrix equal to $\sigma^2I_{M+N}$. In this case, $\mathsf{h}\left[\vec{y}_{0}\right]=0.5\logp{2\pi\e{}\sigma^2}^{N+M}$. For large values of $\sigma$, the distribution of

The next corollary establishes a lower bound on the MSD of primal/dual variables:
\begin{corollary}\label{Coro: FTE-P}
Let $\ES{\norm{\vec{\epsilon}^{\vec{x}}_k}{2}^2}$ and $\ES{\norm{\vec{\epsilon}^{\vec{\lambda}}_k}{2}^2}$ be the MSD of the primal variables and dual variables, respectively, at time $k$ from the optimal solution. Then, we have
\begin{align}
&\log\ES{\norm{\vec{\epsilon}^{\vec{x}}_k}{2}^2}\geq \logp{\frac{\e{1-\frac{1}{M}}}{{2\pi \e{}}}}+ \nonumber\\
&\hspace{0cm}\frac{2}{M}\left(\sum_{i=1}^M\sum_{n=0}^{k-1}{\logp{1+\mu_nU_i^{\min}}}+\mathsf{h}\left[\vec{x}_{0}\right]-\sum_{t=0}^{k-1}\sum_{j=1}^N\log\abs{\mathcal{A}^{\vec{\lambda}}_{j,t}}\right),\nonumber\\
&\log\ES{\norm{\vec{\epsilon}^{\vec{\lambda}}_k}{2}^2}\geq \logp{\frac{\e{1-\frac{1}{N}}}{{2\pi \e{}}}}+ \nonumber\\
&\hspace{2cm}+\frac{2}{N}\paren{\mathsf{h}\left[\vec{\lambda}_{0}\right]-\sum_{t=0}^{k-1}\sum_{i=1}^M\log\abs{\mathcal{A}^{\vec{x}}_{i,t}}},\nonumber
\end{align}
\end{corollary}
\begin{IEEEproof}
The proof is similar to that of Corollary \ref{Coro: FTE} and is omitted to avoid repetition. %lower bound is directly obtained by combining the inequalities \eqref{Eq: MDPLB} and \eqref{Eq: Aux-P} in Appendix \ref{App: DDE-P}.
\end{IEEEproof}

\section{An Optimum Achieving Quantization Scheme}\label{Sec: AOAQ}

\textcolor{black}{In this section, we propose a zoom-in uniform optimum achieving (OA) quantization scheme for the PD algorithm. We refer to this quantization scheme as $\mathcal{Q}_{\rm a}$. We also prove that PD algorithm under the quantization scheme $\mathcal{Q}_{\rm a}$  converges to the optimal solution of the optimization problem \eqref{Prob}. To this end, we assume that the unquantized PD algorithm forms a contraction map with contraction constant $\alpha\in\left[0,\left.1\right)\right.$. We assume that $\alpha$ is known by all agents and NNs.  Let ${\rm Q}_{{\rm a},k}\paren{\cdot}$ and $\delta_k$ denote the quantizer and the quantization step employed by agents and NNs at time $k$. $\delta_k$ is set to $\alpha^{k+1}$. }

\textcolor{black}{At time $k=0$, agent $i$ generates $x^i_0$ according to the uniform distribution on the interval $\left(-L\alpha,L\alpha\right)$ where $L$ is a positive integer. Similarly, NN $j$ generates $\lambda^j_0$ ($1\leq j\leq N$) using the uniform distribution on $\left(-L\alpha,L\alpha\right)$. Next, agents and NNs quantize the initial primal and dual variables, respectively, using the midpoint uniform quantizer on $\left(-L\alpha,L\alpha\right)$ with the quantization step $\delta_0=\alpha$. Thus, the quantizer employed by agents and NNs at time $k=0$, is given by ${\rm Q}_{{\rm a},0}\paren{z}=\left\lfloor \frac{z}{\alpha}\right\rfloor\alpha+\frac{\alpha}{2}$ for $z\in\left(-L\alpha,L\alpha\right)$ where $\left\lfloor \cdot\right\rfloor$ is the floor function. Each agent/NN only needs $\left\lceil \log_2\paren{2L}\right\rceil$ bits to communicate its initial primal/dual variable, respectively, where $\left\lceil \cdot\right\rceil$ is the ceiling function.}

\textcolor{black}{Let $I^{x^i}_{k+1}$ be the interval centered at $C^{x^i}_{k+1}={\rm Q}_{{\rm a},k}\paren{x^i_{k}}+\left\lfloor \frac{x^i_{k+1}-x^i_{k}}{\delta_k}\right\rfloor\delta_k$ with length $2\left\lceil \frac{2}{\alpha}\right\rceil\delta_{k+1}$. It can be shown that $x^i_{k+1}$ belongs to this interval (see the proof of Theorem \ref{Theo: Conv} for more details).  At time $k+1$, a uniform midpoint quantizer is employed to quantize $x^{i}_{k+1}$. To this end, first, agent $i$ transmits $\left\lfloor \frac{x^i_{k+1}-x^i_{k}}{\delta_k}\right\rfloor$ to NNs. Since $\delta_k$ and ${\rm Q}_{{\rm a},k}\paren{x^i_{k}}$ are known by NNs, each NN can compute $C^{x^i}_{k+1}$ using $\left\lfloor \frac{x^i_{k+1}-x^i_{k}}{\delta_k}\right\rfloor$. The knowledge of $C^{x^i}_{k+1}$ allows each NN to update its decoder at time $k+1$. Note that $\left\lfloor \frac{x^i_{k+1}-x^i_{k}}{\delta_k}\right\rfloor$ is an integer which can be transmitted using finite number of bits. Next, agent $i$ computes ${\rm \hat{Q}}_{k+1}\paren{\frac{x^i_{k+1}-C^{x^i}_{k+1}}{\delta_{k+1}}}$ where $x^i_{k+1}\in I^{x^i}_{k+1}$ and ${\rm \hat{Q}}_{k+1}\paren{\cdot}$ is given by 
\begin{eqnarray}\label{Eq: UQant}
\hspace{-2mm}{\rm \hat{Q}}_{k+1}\paren{z}=\!\!\left\{\!\!\!\!
\begin{array}{cc}
\vspace{2mm}\left\lceil \frac{2}{\alpha}\right\rceil-1 &\paren{\left\lceil \frac{2}{\alpha}\right\rceil-1}\leq z\leq \left\lceil \frac{2}{\alpha}\right\rceil\\
\vspace{2mm}\left\lfloor \frac{z}{\delta_{k+1}}\right\rfloor& -\left\lceil \frac{2}{\alpha}\right\rceil\leq z\leq \paren{\left\lceil \frac{2}{\alpha}\right\rceil-1}\\
%-\left\lceil \frac{2}{\alpha}\right\rceil & -\left\lceil \frac{2}{\alpha}\right\rceil\alpha\delta_k\leq z\leq \paren{-\left\lceil \frac{2}{\alpha}\right\rceil+1}\alpha\delta_k
\end{array}
\right.
\end{eqnarray}
Finally, agent $i$ transmits ${\rm \hat{Q}}_{k+1}\paren{\frac{x^i_{k+1}-C^{x^i}_{k+1}}{\delta_{k+1}}}$ to NNs using $\left\lceil \log_2\paren{2\left\lceil \frac{2}{\alpha}\right\rceil}\right\rceil$ bits. Then, each NN constructs the quantized version of $x^i_{k+1}$ using ${\rm Q}_{{\rm a},k+1}\paren{x^i_{k+1}}=C^{x^i}_{k+1}+{\rm \hat{Q}}_{k+1}\paren{x^i_{k+1}-C^{x^i}_{k+1}}\delta_{k+1}+\frac{\delta_{k+1}}{2}$.}

\textcolor{black}{Consider the interval $I^{\lambda^j}_{k+1}$ which is centered at $C^{\lambda^j}_{k+1}={\rm Q}_{{\rm a},k}\paren{\lambda^j_{k}}+\left\lfloor \frac{\lambda^j_{k+1}-\lambda^j_{k}}{\delta_k}\right\rfloor\delta_k$ with the length  $2\left\lceil \frac{2}{\alpha}\right\rceil\delta_{k+1}$. It can be shown that the $\lambda^j_{k+1}$ belongs to $I^{\lambda^j}_{k+1}$ (see the proof of Theorem \ref{Theo: Conv} for more details). At time $k+1$, NN $j$, first, broadcasts $\left\lfloor \frac{\lambda^j_{k+1}-\lambda^j_{k}}{\delta_k}\right\rfloor$ to agents which allows the agents to construct $C^{\lambda^j}_{k+1}$. Then, it broadcasts  ${\rm \hat{Q}}_{k+1}\paren{\frac{\lambda^j_{k+1}-C^{\lambda^j}_{k+1}}{\delta_{k+1}}}$ to all agents where ${\rm \hat{Q}}_{k+1}\paren{\cdot}$ is given by \eqref{Eq: UQant}. Finally, agents construct the quantized version of $\lambda^j_{k+1}$ as ${\rm Q}_{{\rm a},k+1}\paren{\lambda^j_{k+1}}=C^{\lambda^j}_{k+1}+{\rm \hat{Q}}_{k+1}\paren{\frac{\lambda^j_{k+1}-C^{\lambda^j}_{k+1}}{\delta_{k+1}}}\delta_{k+1}+\frac{\delta_{k+1}}{2}$.}

\textcolor{black}{The next theorem shows that the quantized PD algorithm under $\mathcal{Q}_{\rm a}$ converges to the optimal solution.
\begin{theorem}\label{Theo: Conv}
The PD algorithm under the quantization scheme $\mathcal{Q}_{\rm a}$ converges  exponentially to the optimal solution of the optimization problem \eqref{Prob}. 
\end{theorem}
\begin{IEEEproof}
Please see Appendix \ref{App:Conv}.
\end{IEEEproof}}
\section{Numerical results}\label{Sec: NR}
\textcolor{black}{This section presents our numerical results illustrating the behavior of mean square distance (MSD) of primal-dual (PD) variables from the optimal solution, with time index $k$ for a quadratic network utility maximization (NUM) problem with 10 agents and 5 constraints under the quantization scheme $\mathcal{Q}_{a}$. %We assume that the objective function of agent $i$ is given by $U_i\paren{x^i}=\alpha_i\paren{x^i\mathsf{\tan}^{-1}\paren{x^i}-0.5\logp{1+{x^i}^2}}$ where $\alpha_i<0$. It can be easily verified that the objective function of agent $i$ is a concave function of $x^i$. 
%Let $\mathcal{Q}_{\rm u}$ denote the uniform quantization scheme  employed by agents and the network manager (NM) to quantize the primal and dual variables, respectively. Under the quantization scheme $\mathcal{Q}_{\rm u}$, at time $k$, agent $i$ will select an interval of length $d_k$ containing $x^i_k$. Let $I_{k}^{x^i}\in \R$ be such an interval. Then, agent $i$ quantizes $I_{k}^{x^i}$ in to $2^R$ equal length sub-intervals and transmits the representative of the sub-interval containing $x^i_k$ to the NM as its decision variable at time $k$. Here, $R$ is the communication data rate allocated for transmitting $x^i$ to the NM.
%Similarly, at time $k$, the NM selects the interval $I^{\lambda^j}_k$ with length $d_k$ which contains $\lambda^j_k$, and divides it into $2^R$ equal length sub-intervals where $R$ is the data rate allocated for transmitting $\lambda^j$ to agents. Then, the NM transmits the representative of the sub-interval containing $\lambda^j_k$ to all agents as the value of the $j$th dual variable at time $k$. We assume that $d_k$ decays with time according to $d_k=\beta_T d_{k-1}$ where $\beta_T\in\left(0,1\right)$.
In our numerical results, the step-size of the PD algorithm, $\mu_k$, is set to $0.019$ for all $k$, $\alpha=0.9495$, $L=5$. The initial PD variables are independently drawn according to the uniform distribution on $\left(-L\alpha,L\alpha\right)$. For $k\geq 1$,  3 bits is used to quantize each primal/dual variable.}

\begin{figure}[t]
\centering
\subfigure[]
{\includegraphics[scale=0.6]{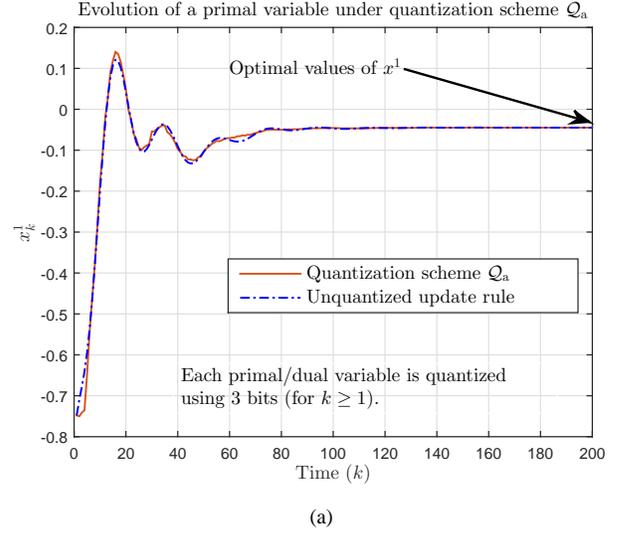}
\label{F1}}
%\subfigure[]
%{\includegraphics[scale=0.6]{Dual_time.eps}
%\label{F2}}
\caption{\small Time evolution of the primal variable $x^1$ under the quantization scheme $\mathcal{Q}_{\rm a}$ and unquantized PD update rule. The initial PD variables are the same in both graphs.}
\label{F1-2}
\end{figure}
\textcolor{black}{Fig. \ref{F1-2} shows the time evolution of the primal variable $x^1$ under the quantization scheme $\mathcal{Q}_{\rm a}$ and unquantized PD update rule. The initial PD variables are the same for both graphs. According to Fig. \ref{F1-2}, the trajectories of $x^1$, under both quantization scheme $\mathcal{Q}_{\rm a}$ and the unquantized PD update rule, converge to the optimal value of ${x^1}$ as the time index $k$ becomes large. The same behavior continues to hold for other primal/dual variables. }

\begin{figure}[!t]
\centering{\includegraphics[scale=0.6]{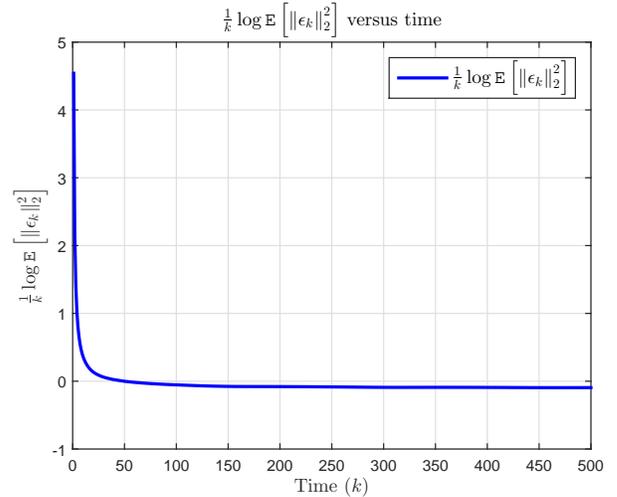}}
\caption{$\log$-MSD divided by $k$ for primal-dual variables versus time index ($k$) under the quantization scheme $\mathcal{Q}_{\rm a}$.} \label{F3}
\end{figure}
\textcolor{black}{Fig. \ref{F3} depicts $\log$-MSD divided by $k$ for the PD variables , \emph{i.e.,} $\frac{1}{k}\log\ES{\norm{\vec{\epsilon}_k}{2}^2}$, as a function of time index $k$. The lower bounds on the distance decay exponent (DDE) of PD variables, predicted by Theorem \ref{Theo: DDE} and Theorem \ref{Theo: EDE-New},  are equal to $ -6.24$ and $-4.64$, respectively. As Fig. \ref{F3} shows, $\frac{1}{k}\log\ES{\norm{\vec{\epsilon}_k}{2}^2}$ stays above the predicted values by Theorems \ref{Theo: DDE} and \ref{Theo: EDE-New} as $k$ becomes large.}
\section{Conclusions}\label{Sec: Conc}
In this paper, we have studied the convergence behavior of the  quantized primal-dual (PD) algorithm in solving network utility maximization problems. First, using the information-theoretic notion of entropy power, we established universal bounds on the fastest speed of exponential mean square convergence of PD, primal and dual variables to the optimal solution under optimum achieving quantization schemes. Here, our results provide universal trade-offs between the speed of convergence of the quantized PD algorithm, data rate under the quantization, objective functions of agents,  the number of agents and the number of constraints.   Next, we established universal lower bounds on the mean square distance of PD, primal and dual variables from the optimal solution of the NUM problem for any finite time index. 
%%%%%%%%%%%%%%%%%%%%%%%%%%%%%%%%%%%%%%%%%%%%%%%%%%%%%%%%%%%%%%%%%%%%%%%%%%%%%%%%%%%%%%%%%%%%%%%%%%%%%%%%%%%%%%%%%%%%%%%%%%%%%%%%%
\appendices

\section{Proof of Theorem \ref{Theo: DDE}}\label{App: DDE}

This appendix presents the main steps of the proof of Theorem \ref{Theo: DDE}. To this end, first, the notion of conditional differential entropy power of a random vector is defined. Then, we use the notion of entropy power to establish a universal lower bound on the DDE of the PD variables. The differential entropy power of the random vector $\vec{z}\in \R^{N+M} $ conditioned on the event $A=a$, denoted by $\CEP{\vec{z}}{A=a}$, is defined as
\begin{eqnarray}
\CEP{\vec{z}}{A=a}=\frac{1}{{2\pi \e{}}}\e{\frac{2}{M+N}\CENT{\vec{z}}{A=a}},\nonumber
\end{eqnarray}
where $\CENT{\vec{z}}{A=a}$ is the conditional differential entropy of $\vec{z}$ given  $A=a$ defined as
\begin{eqnarray}
\CENT{\vec{z}}{A=a}=-\int\logp{p\paren{\left.\vec{z}\right|A=a}}p\paren{\left.\vec{z}\right|A=a}d\vec{z},\nonumber
\end{eqnarray}
where $p\paren{\left.\vec{z}\right|A=a}$ is the conditional distribution of $\vec{z}$ given $A=a$.
Using the entropy maximizing property of Gaussian distributions, the conditional entropy power of $\vec{z}$ given $A=a$ can be upper bounded \cite{NE04} as
\begin{eqnarray}\label{Eq: CEP-UB}
\CEP{\vec{z}}{A=a}\leq \e{1/\paren{M+N}-1}\CES{\norm{\vec{z}}{2}^2}{A=a},
\end{eqnarray}
where $\CES{\vec{z}}{A=a}$ is conditional expectation of $\vec{z}$ given $A=a$. Let $\mathsf{E}_A\left[\CEP{\vec{z}}{A=a}\right]$ denote the average conditional entropy power of $\vec{z}$ given $A=a$. Using \eqref{Eq: CEP-UB}, $\mathsf{E}_A\left[\CEP{\vec{z}}{A=a}\right]$ can be upper bounded as
\begin{eqnarray}\label{Eq: ACEP-UB}
{\mathsf E}_{A}\left[\CEP{\vec{z}}{A}\right]\leq \e{1/\paren{M+N}-1}\ES{\norm{\vec{z}}{2}^2}.
\end{eqnarray}

Next, the inequality \eqref{Eq: ACEP-UB} is used to establish the universal lower bound on the DDE of the PD variables under OA quantization schemes. \textcolor{black}{To this end, let $\mathcal{D}_{k-1}=\left\{\hat{Q}_n=\vec{\hat{q}}_n\right\}_{n=0}^{k-1}$ where $\hat{Q}_n=\left[\hat{Q}^{\vec{x}}_{1,n},\cdots,\hat{Q}^{\vec{x}}_{M,n},\hat{Q}^{\vec{\lambda}}_{1,n},\cdots,\hat{Q}^{\vec{\lambda}}_{N,n}\right]$ and $\vec{\hat{q}}_n$ is a possible realization of $\hat{Q}_n$.} Using \eqref{Eq: ACEP-UB}, $\ES{\norm{\vec{\epsilon}_k}{2}^2}$ can be lower bounded as
\begin{eqnarray}\label{Eq: MS-LB}
\ES{\norm{\vec{\epsilon}_k}{2}^2}&\geq &\e{1-\frac{1}{M+N}}\ES{\CEP{\vec{\epsilon}_k}{\mathcal{D}_{k-1}}}\nonumber\\
&\stackrel{\paren{*}}{\geq}&\frac{\e{1-\frac{1}{M+N}}}{{2\pi \e{}}}\e{\frac{2}{M+N}\ES{\CENT{\vec{\epsilon}_k}{\mathcal{D}_{k-1}}}},
\end{eqnarray}
where $\paren{*}$ is obtained using the Jensen inequality. The term $\CENT{\vec{\epsilon}_k}{\mathcal{D}_{k-1}}$ on the right hand side of \eqref{Eq: MS-LB} can be expanded as
\begin{align}\label{Eq: TI}
\CENT{\vec{\epsilon}_k}{\mathcal{D}_{k-1}}&=\CENT{\vec{y}_k-\vec{y}^\star}{\mathcal{D}_{k-1}}\nonumber\\
&\stackrel{\paren{*}}{=}\CENT{\vec{y}_k}{\mathcal{D}_{k-1}},
\end{align}
where $\paren{*}$ follows from the translation invariance property of differential entropy \textcolor{black}{as $\vec{y}^\star$ is a constant vector (see \cite{Cover} Theorem 8.6.3 page 253)}. 

The next lemma establishes a useful expression between $\CENT{\vec{y}_n}{\mathcal{D}_{k-1}}$ and $\CENT{\vec{y}_{n-1}}{\mathcal{D}_{k-1}}$ for $n\leq k$, which is used to further expand $\CENT{\vec{y}_k}{\mathcal{D}_{k-1}}$.

\begin{lemma}\label{Lem: Entropy-Equality}
For $n\leq k$, $\CENT{\vec{y}_n}{\mathcal{D}_{k-1}}$ can be expanded as
\begin{align}
&\CENT{\vec{y}_{n}}{\mathcal{D}_{k-1}}=\CENT{\vec{y}_{n-1}}{\mathcal{D}_{k-1}}+\nonumber\\
&\mathsf{E}\left[\sum_{j=1}^M\log\left(1+\mu_{n-1}\left.\frac{d ^2}{d {x^j}^2}U_{j}\left(\vec{x}^{j}_{n-1}\!\right)\right)\right|\mathcal{D}_{k-1}\right]
\end{align}
\end{lemma}
\begin{IEEEproof}
Let $\tilde{x}^i_n={x}^{i}_{n}+\mu_{n} \paren{\frac{d}{d {x}^{i}}U_i\paren{{x}^{i}_{n}}}$ and $\vec{\tilde{x}}_n=\left[\tilde{x}^i_1,\cdots,\tilde{x}^i_M\right]^\top$. Let $\vec{\tilde{y}}_n$ be the vector concatenation of $\vec{\tilde{x}}_n$ and $\vec{\lambda}_n$. \textcolor{black}{This lemma is proved in two steps. First, it is shown that the conditional differential entropy of $\vec{y}_n$ given $\mathcal{D}_k$ is equal to that of $\vec{\tilde{y}}_{n-1}$ given $\mathcal{D}_k$ (see \eqref{Eq: ENT-Eq-1}). Next, a relation between the conditional differential entropy of $\vec{\tilde{y}}_{n-1}$ given $\mathcal{D}_k$ and that of $\vec{y}_{n-1}$ given $\mathcal{D}_k$ is established. Note that, $\CENT{\vec{y}_n}{\mathcal{D}_{k-1}}$ can be written as
\begin{eqnarray}\label{Eq: ENT-Eq-1}
\CENT{\vec{y}_n}{\mathcal{D}_{k-1}}&=&\CENT{\vec{x}_n,\vec{\lambda}_n}{\mathcal{D}_{k-1}}\nonumber\\
&\stackrel{*}{=}&\CENT{\vec{\tilde{x}}_{n-1},\vec{\lambda}_{n-1}}{\mathcal{D}_{k-1}}\nonumber\\
&=&\CENT{\vec{\tilde{y}}_{n-1}}{\mathcal{D}_{k-1}}
\end{eqnarray}
where $\paren{*}$ follows from the translation invariance property of the differential entropy and the fact that  $Q_{k-1}$ is fixed given $\mathcal{D}_{k-1}=\left\{\hat{Q}_n=\vec{\hat{q}}_n\right\}_{n=0}^{k-1}$.} Next, we derive an expression for the probability density function (PDF) of $\vec{\tilde{y}}_n$ in terms of the PDF of $\vec{y}_{n}$. Let $p_{\vec{\tilde{y}}_{n}}\paren{\vec{y}\left|\mathcal{D}_{k-1}\right.}$ and $p_{\vec{y}_{n}}\paren{\vec{y}\left|\mathcal{D}_{k-1}\right.}$ to denote the PDFs of $\vec{\tilde{y}}_n$ and $\vec{y}_{n}$, respectively, conditioned on $\mathcal{D}_{k-1}$.  Let $\vec{F}\paren{\cdot}$ represent the mapping between $\vec{\tilde{y}}_n$ and $\vec{y}_n$, \emph{i.e.,} $\vec{\tilde{y}}_n=\vec{F}\paren{\vec{y}_n}$.  Note that $0<1+\mu_{n}\frac{d^2}{d {x^i}^2}U_i\paren{x^i}<1$ since  $0<\mu_{n}< \min_i\frac{1}{\abs{U^{\rm min}_i}}$ which implies that the mapping $\vec{F}\paren{\cdot}$ is invertible. 
\textcolor{black}{Thus, the change-of-variables formula for invertible diffeomorphisms of random vectors (see {\em e.g.}, (4.63) in \cite{leon-garcia}) can be applied to write}

	\begin{eqnarray}\label{Eq: PDF-R}
	p_{\vec{\tilde{y}}_{n-1}}\!\!\paren{\vec{y}\left|\mathcal{D}_{k-1}\right.}=\frac{1}{\det J_{\vec{F}}\!\!\left[\vec{F}^{-1}\!\!\paren{\vec{y}}\right]}p_{\vec{y}_{n-1}}\!\!\paren{\vec{F}^{-1}\paren{\vec{y}}\left|\mathcal{D}_{k-1}\right.},\nonumber\\
	%&=&\left.\prod_{j=1}^M\frac{1}{{1+\mu_N\frac{d^2}{d {x^i}^2}U_i\paren{x_i,\vec{a}^{-i}_{k-1}}}}\right|_{F^{-1}_i\paren{x^i}}p_{\act{\vec{x}}{}{k-1}}\paren{F^{-1}\paren{\vec{x}}\left|\left\{\gquant{\act{\vec{x}}{}{i}}{i}=\vec{a}_i\right\}_{i=0}^{N-1}\right.}
	\end{eqnarray}
	where $J_{\vec{F}}\left[\vec{x}\right]$ is Jacobian of  $\vec{F}\paren{\vec{x}}$ evaluated at $\vec{x}$. Using \eqref{Eq: PDF-R}, the conditional entropy of $\vec{\tilde{y}}_{n-1}$ given $\mathcal{D}_{k-1}$ can be written as
\begin{align}
&\CENT{\vec{\tilde{y}}_{n-1}}{\mathcal{D}_{k-1}}\nonumber\\
&=\int\logp{\det J_{\vec{F}}\left[\vec{F}^{-1}\paren{\vec{y}}\right]}\frac{1}{\det J_{\vec{F}}\left[\vec{F}^{-1}\paren{\vec{y}}\right]}\nonumber\\
&\hspace{5cm}p_{\vec{y}_{n-1}}\paren{F^{-1}\paren{\vec{y}}\left|\mathcal{D}_{k-1}\right.}d\vec{y}\nonumber\\
&-\int\logp{p_{\vec{y}_{n-1}}\paren{F^{-1}\paren{\vec{y}}\left|\mathcal{D}_{k-1}\right.}}\frac{1}{\det J_{\vec{F}}\left[\vec{F}^{-1}\paren{\vec{y}}\right]}\nonumber\\
&\hspace{5cm}p_{\vec{y}_{n-1}}\paren{\vec{F}^{-1}\paren{\vec{y}}\left|\mathcal{D}_{k-1}\right.}d\vec{y},\nonumber\\
&\stackrel{\paren{*}}{=}\int\logp{\det J_{\vec{F}}\left[\vec{z}\right]}p_{\vec{y}_{n-1}}\paren{\vec{z}\left|\mathcal{D}_{k-1}\right.}d\vec{z}\nonumber\\
&\hspace{1cm}-\int\logp{p_{\vec{y}_{n-1}}\paren{\vec{z}\left|\mathcal{D}_{k-1}\right.}}p_{\vec{y}_{n-1}}\paren{\vec{z}\left|\mathcal{D}_{k-1}\right.}d\vec{z},\nonumber\\
&=\!\!\sum_{j=1}^M\mathsf{E}\!\left[\left.\!\log\!\!\left(\!1\!+\!\mu_{n-1}\frac{d ^2}{d {x^j}^2}U_{j}\!\!\paren{x^j_{n-1}}\right)\right|{\mathcal{D}_{k-1}}\right]\nonumber\\
&\hspace{1cm}+\CENT{\vec{y}_{n-1}}{\mathcal{D}_{k-1}},
\end{align}	
where $\paren{*}$ follows from the change of variable $\vec{z}=F^{-1}\paren{\vec{x}}$.
\end{IEEEproof}
Using Lemma \ref{Lem: Entropy-Equality}, $\CENT{\vec{y}_{k}}{\mathcal{D}_{k-1}}$ can be further expanded as
\begin{align}\label{Eq: Aux-AE-LB}
&\CENT{\vec{y}_{k}}{\mathcal{D}_{k-1}}=\CENT{\vec{y}_{0}}{\mathcal{D}_{k-1}}+\nonumber\\
&\sum_{j=1}^M\sum_{n=0}^{k-1}\CES{\logp{1+\mu_n\frac{d ^2}{d {x^j}^2}U_{j}\paren{x^{j}_{n}}}}{\mathcal{D}_{k-1}}
\end{align}
Using \eqref{Eq: Aux-AE-LB}, $\ES{\CENT{\vec{y}_{k}}{\mathcal{D}_{k-1}}}$ can be written as
\begin{align}\label{Eq: Aux-0}
&\ES{\CENT{\vec{y}_{k}}{\mathcal{D}_{k-1}}}\nonumber\\
&\hspace{0cm}=\sum_{j=1}^M\sum_{n=0}^{k-1}\ES{\logp{1+\mu_n\frac{d ^2}{d {x^j}^2}U_{j}\paren{x^{j}_{n}}}}+\ES{\CENT{\vec{y}_{0}}{\mathcal{D}_{k-1}}},
\end{align}
The following lemma, adapted from  \cite{NE04},
establishes a lower bound on $\ES{\CENT{\vec{y}_{k}}{\mathcal{D}_{k-1}}}$:
\begin{lemma}\label{Lem: NE}
The average conditional entropy of $\vec{y}_0$ given $\mathcal{D}_{k-1}$, \emph{i.e.,} $\ES{\CENT{\vec{y}_{0}}{\mathcal{D}_{k-1}}}$, can be lower bounded as
\begin{eqnarray}
\ES{\CENT{\vec{y}_{0}}{\mathcal{D}_{k-1}}}\!\geq \!\mathsf{h}\left[\vec{y}_{0}\right]\!-\!\sum_{t=0}^{k-1}\!\paren{\!\!\!\paren{\sum_{i=1}^M\log\abs{\mathcal{A}^{\vec{x}}_{i,t}}}\!\!+\!\!\sum_{j=1}^N\log\abs{\mathcal{A}^{\vec{\lambda}}_{j,t}}\!\!}.\nonumber
\end{eqnarray}
\end{lemma}
\begin{IEEEproof}
 Follows directly from the first inequality in appendix C in  \cite{NE04}; alternatively, it can be derived from (8.48) and (8.89) in \cite{Cover}.
\end{IEEEproof}

Applying Lemma \ref{Lem: NE} to \eqref{Eq: Aux-0}, we have
\begin{align}\label{Eq: Aux-1}
&\ES{\CENT{\vec{y}_{k}}{\mathcal{D}_{k-1}}}\nonumber\\
&\hspace{1cm}\geq \sum_{j=1}^M\sum_{n=0}^{k-1}\ES{\logp{1+\mu_n\frac{d ^2}{d {x^j}^2}U_{j}\paren{x^{j}_{n}}}}\nonumber\\
&\hspace{1cm}+\mathsf{h}\left[\vec{y}_{0}\right]-\sum_{t=0}^{k-1}\paren{\paren{\sum_{i=1}^M\log\abs{\mathcal{A}^{\vec{x}}_{i,t}}}+\sum_{j=1}^N\log\abs{\mathcal{A}^{\vec{\lambda}}_{j,t}}},
\end{align}
Since $\vec{x}_0$ and $\vec{\lambda}_0$ are independent, the differential entropy of $\vec{y}_0$ can be written as $\mathsf{h}\left[\vec{y}_{0}\right]=\mathsf{h}\left[\vec{x}_{0}\right]+\mathsf{h}\left[\vec{\lambda}_{0}\right]$ which implies that $\vec{y}_0$ has finite differential entropy. Using \eqref{Eq: MS-LB}, \eqref{Eq: TI}, \eqref{Eq: Aux-1} and the fact that $\vec{y}_{0}$ has a finite entropy,  the DDE can be lower bounded as
\begin{align}\label{Eq: Aux-2}
&\liminf_{k\longrightarrow\infty}\frac{1}{k}\log\ES{\norm{\vec{\epsilon}_k}{2}^2}\geq \frac{2}{M+N} \paren{\liminf_{k\longrightarrow \infty}\sum_{j=1}^M\frac{1}{k}\sum_{n=0}^{k-1}\right.\nonumber\\
&\left.\ES{\logp{1+\mu_n\frac{d ^2}{d {x^j}^2}U_{j}\paren{x^{j}_{n}}}}-R_{\mathcal{Q}}}.
\end{align}
In the next lemma, we study the asymptotic behavior of the first term in the right hand side of equation \eqref{Eq: Aux-2}.
\begin{lemma}\label{Lem: Limit}
Consider the primal-dual update rule \eqref{EQ: QU} under an OA quantization scheme. Then, we have
\begin{align}
\lim_{k\longrightarrow \infty}\sum_{j=1}^M&\frac{1}{k}\sum_{n=0}^{k-1}\ES{\logp{1+\mu_n\frac{d ^2}{d {x^j}^2}U_{j}\paren{{x}^{j}_{n}}}}=\nonumber\\
&\sum_{j=1}^M\logp{1+\mu^\star\frac{d ^2}{d {x^j}^2}U_{j}\paren{{x^{j}}^\star}}.\nonumber
\end{align}
\end{lemma}
\begin{IEEEproof}
To prove this lemma, first we define the events $E$ and $E_n$ as
\begin{eqnarray}
E&=&\left\{\norm{\vec{x}-\vec{x}^\star}{2}\leq \delta, \abs{\mu-\mu^\star}\leq \delta\right\},\nonumber\\
E_n&=&\left\{\norm{\vec{x}_{n}-\vec{x}^\star}{2}\leq \delta, \abs{\mu_n-\mu^\star}\leq \delta\right\}.\nonumber
\end{eqnarray}
respectively, where $\delta>0$. Then, we have
\begin{align}\label{Eq: Expansion}
&\frac{1}{k}\sum_{n=0}^{k-1}\sum_{j=1}^M\ES{\logp{1+\mu_n\frac{d ^2}{d {x^j}^2}U_{j}\paren{{x}^{j}_{n}}}}\nonumber\\
&=\frac{1}{k}\sum_{n=0}^{k-1}\sum_{j=1}^M\ES{\logp{1+\mu_n\frac{d ^2}{d {x^j}^2}U_{j}\paren{{x}^{j}_{n}}}\paren{\I{E_n}+\I{E^{c}_n}}},\nonumber\\
&\geq \inf_{\left\{\vec{x},\mu\right\}\in E}\sum_{j=1}^M\logp{1+\mu \frac{d ^2}{d {x^j}^2}U_j\paren{{x}^{j}}}\ES{\frac{1}{k}\sum_{n=0}^{k-1}\I{E_n}}\nonumber\\
&\hspace{2cm}+\sum_{j=1}^M\logp{1+\inf_n\mu_nU_j^{\rm min}}\ES{\frac{1}{k}\sum_{n=0}^{k-1}\I{E^{c}_n}},
\end{align}
where $E^{c}_n$ is the complement of the event $E_n$. Recall that the quantization scheme $\mathcal{Q}$ is an OA quantization scheme. Therefore, we have $\lim_{k\longrightarrow \infty} {\vec{x}}_{k}=\vec{x}^\star$ for any initial vector $\vec{x}_0$ in the support set of $p_{\vec{x}_0}\paren{\vec{x}}$ which implies $\I{E_n}\longrightarrow 1$ almost surely and $\frac{1}{k}\sum_{n=0}^{k-1}\I{E_n}\longrightarrow 1$ almost surely. Applying Fatou's Lemma to \eqref{Eq: Expansion}, we have $ \liminf_{k\longrightarrow\infty}\ES{\frac{1}{k}\sum_{n=0}^{k-1}\I{E_n}}\geq 1$. Using Lebesgue dominated convergence Theorem and the fact that $\frac{1}{k}\sum_{n=0}^{k-1}\I{E^{c}_n}\longrightarrow 0$ almost surely, we have $\lim_{k\longrightarrow \infty}\ES{\frac{1}{k}\sum_{n=0}^{k-1}\I{E^{c}_n}}=0$. Hence,
\begin{align}
&\hspace{-1cm}\liminf_{k\longrightarrow \infty}\frac{1}{k}\sum_{n=0}^{k-1}\sum_{j=1}^M\ES{\logp{1+\mu_n\frac{d ^2}{d {x^j}^2}U_{j}\paren{{x}^{j}_{n}}}}
\geq\nonumber\\
& \hspace{1cm}\inf_{\left\{\vec{x},\mu\right\}\in E}\sum_{j=1}^M\logp{1+\mu \frac{d^2}{d {x^j}^2} U_j\paren{\vec{x}^{j}}}.\nonumber
\end{align}
Since, $\delta>0$ is arbitrary, we have
\begin{align}
\liminf_{k\longrightarrow \infty}\frac{1}{k}&\sum_{n=0}^{k-1}\sum_{j=1}^M\ES{\logp{1+\mu_n\frac{d ^2}{d {x^j}^2}U_{j}\paren{{x}^{j}_{n}}}}
\geq\nonumber\\
& \sum_{j=1}^M\logp{1+\mu^\star \frac{d^2}{d {x^j}^2} U_j\paren{{{x}^{j}}^\star}}\nonumber
\end{align}
To prove the other direction, note that
\begin{align}
&\frac{1}{k}\sum_{n=0}^{k-1}\sum_{j=1}^M\ES{\logp{1+\mu_n\frac{d ^2}{d {x^j}^2}U_{j}\paren{{x}^{j}_{n}}}}\leq\nonumber\\
& \sup_{\left\{\vec{x},\mu\right\}\in E}\sum_{j=1}^M\logp{1+\mu \frac{d^2}{d {x^j}^2} U_j\paren{{x}^{j}}}\ES{\frac{1}{k}\sum_{n=0}^{k-1}\I{E_n}}\nonumber\\
&+\sum_{j=1}^M\logp{1+\sup_n\mu_nU^{\rm max}_j}\ES{\frac{1}{k}\sum_{n=0}^{k-1}\I{E^{c}_n}}.\nonumber
\end{align}
Following similar steps as before, it can be easily shown that
\begin{align}
\limsup_{k\longrightarrow \infty}&\frac{1}{k}\sum_{n=0}^{k-1}\sum_{j=1}^M\ES{\logp{1+\mu_n\frac{d ^2}{d {x^j}^2}U_{j}\paren{{x}^{j}_{n}}}}
\leq\nonumber\\
& \sum_{j=1}^M\logp{1+\mu^\star\frac{d^2}{d {x^j}^2} U_j\paren{{{x}^{j}}^\star}},\nonumber
\end{align}
which completes the proof.
\end{IEEEproof}
Applying Lemma \ref{Lem: Limit} to \eqref{Eq: Aux-2}, we have
\begin{align}
&\liminf_{k\rightarrow\infty}\frac{1}{k}\log\ES{\norm{\vec{\epsilon}_k}{2}^2}\nonumber\\
&\hspace{1cm}\geq\frac{2}{N+M}\paren{\sum_{i=1}^m\log\paren{1+\mu^\star \frac{d^2}{d {x^i}^2}U_i\paren{{x^i}^\star}}-R_{\mathcal{Q}}}.
\end{align}
which completes the proof.
%%%%%%%%%%%%%%%%%%%%%%%%%%%%%%%%%%%%%%%%%%%%%%%%%%%%%%%%%%%%%%%%%%%%%%%%%%%%%%%%%%%%%%%%%%%%%%%%%%%%%%%%%%%%%%%%%%%%%%%%%%%%%%%%
\section{Proof of Theorem \ref{Theo: DDE-P}}\label{App: DDE-P}
The proof essentially follows from similar steps as the proof of Theorem \ref{Theo: DDE}. Here, we state the main steps of the proof for the sake of completeness. Let $\mathcal{D}^{\vec{\lambda}}_{k-1}=\left\{\hat{Q}^{\vec{\lambda}}_{n}=\vec{\hat{q}}^{\vec{\lambda}}_n\right\}_{n=0}^{k-1}$ where  $\hat{Q}^{\vec{\lambda}}_n=\left[\hat{Q}^{\vec{\lambda}}_{1,n},\cdots,\hat{Q}^{\vec{\lambda}}_{N,n}\right]$ and  $\vec{\hat{q}}^{\vec{\lambda}}_n$ is a realization of $\hat{Q}^{\vec{\lambda}}_n$. Then using \eqref{Eq: ACEP-UB}, $\ES{\norm{\vec{\epsilon}_k^{\vec{x}}}{2}^2}$ can be lower bounded as
\begin{eqnarray}\label{Eq: MDPLB}
\ES{\norm{\vec{\epsilon}_k^{\vec{x}}}{2}^2}&\geq &\frac{\e{1-\frac{1}{M}}}{{2\pi \e{}}}\e{\frac{2}{M}\ES{\CENT{\vec{\epsilon}_k^{\vec{x}}}{\mathcal{D}^{\vec{\lambda}}_{k-1}}}}\nonumber\\
&=&\frac{\e{1-\frac{1}{M}}}{{2\pi \e{}}}\e{\frac{2}{M}\ES{\CENT{\vec{x}_k}{\mathcal{D}^{\vec{\lambda}}_{k-1}}}},
\end{eqnarray}
where the equality follows from translation invariance of the differential entropy. Using similar technique as the proof of Lemma \ref{Lem: Entropy-Equality} in Appendix \ref{App: DDE}, one can show
\begin{align}
&\CENT{\vec{x}_{n}}{\mathcal{D}^{\vec{\lambda}}_{k-1}}=\CENT{\vec{x}_{n-1}}{\mathcal{D}^{\vec{\lambda}}_{k-1}}+\nonumber\\
&\mathsf{E}\left[\sum_{j=1}^M\log\left(1+\mu_{n-1}\left.\frac{d ^2}{d {x^j}^2}U_{j}\left(\vec{x}^{j}_{n-1}\!\right)\right)\right|\mathcal{D}^{\vec{\lambda}}_{k-1}\right]
\end{align}
for $n\leq k$. Thus, we have
\begin{align}
&\CENT{\vec{x}_{k}}{\mathcal{D}^{\vec{\lambda}}_{k-1}}=\CENT{\vec{x}_{0}}{\mathcal{D}^{\vec{\lambda}}_{k-1}}+\nonumber\\
&\sum_{j=1}^M\sum_{n=0}^{k-1}\CES{\logp{1+\mu_n\frac{d ^2}{d {x^j}^2}U_{j}\paren{x^{j}_{n}}}}{\mathcal{D}^{\vec{\lambda}}_{k-1}}
\end{align}
and
\begin{align}
&\ES{\CENT{\vec{x}_{k}}{\mathcal{D}^{\vec{\lambda}}_{k-1}}}\nonumber\\
&\hspace{1cm}=\sum_{j=1}^M\sum_{n=0}^{k-1}\ES{\logp{1+\mu_n\frac{d ^2}{d {x^j}^2}U_{j}\paren{x^{j}_{n}}}}\nonumber\\
&\hspace{2cm}+\ES{\CENT{\vec{x}_{0}}{\mathcal{D}^{\vec{\lambda}}_{k-1}}},
\end{align}
Using Lemma \ref{Lem: NE} from Appendix \ref{App: DDE}, $\ES{\CENT{\vec{x}_{k}}{\mathcal{D}^{\vec{\lambda}}_{k-1}}}$ can be lower bounded as
\begin{align}\label{Eq: Aux-P}
&\ES{\CENT{\vec{x}_{k}}{\mathcal{D}^{\vec{\lambda}}_{k-1}}}\nonumber\\
&\hspace{1cm}\stackrel{\paren{*}}{\geq} \sum_{j=1}^M\sum_{n=0}^{k-1}\ES{\logp{1+\mu_n\frac{d ^2}{d {x^j}^2}U_{j}\paren{x^{j}_{n}}}}\nonumber\\
&\hspace{2cm}+\mathsf{h}\left[\vec{x}_{0}\right]-\sum_{t=0}^{k-1}\sum_{j=1}^N\paren{\log\abs{\mathcal{A}^{\vec{\lambda}}_{j,t}}},
\end{align}
Thus, we have
\begin{align}
&\liminf_{k\longrightarrow\infty}\frac{1}{k}\log\ES{\norm{\vec{\epsilon}^{\vec{x}}_k}{2}^2}\geq \frac{2}{M} \paren{\liminf_{k\longrightarrow \infty}\sum_{j=1}^M\frac{1}{k}\sum_{n=0}^{k-1}\right.\nonumber\\
&\left.\ES{\logp{1+\mu_n\frac{d ^2}{d {x^j}^2}U_{j}\paren{x^{j}_{n}}}}-R_{\vec{\lambda}}}.
\end{align}
The proof is complete by appealing to Lemma \ref{Lem: Limit} in Appendix \ref{App: DDE}.
%%%%%%%%%%%%%%%%%%%%%%%%%%%%%%%%%%%%%%%%%%%%%%%%%%%%%%%%%%%%%%%%%%%%%%%%%%%%%%%%%%%%%%%%%%%%%%%%%%%%%%%%%%%%%%%%%%%%%%%%%%%%%%%%
\section{Proof of Theorem \ref{Theo: DDE-D}}\label{App: DDE-D}
Here, we state the main steps of the proof of Theorem \ref{Theo: DDE-D}. \textcolor{black}{Let $\mathcal{D}^{\vec{x}}_{k-1}=\left\{\hat{Q}^{\vec{x}}_{n}=\vec{\hat{q}}^{\vec{x}}_n\right\}_{n=0}^{k-1}$ where $\hat{Q}^{\vec{x}}_n=\left[\hat{Q}^{\vec{x}}_{1,n},\cdots,\hat{Q}^{\vec{x}}_{M,n}\right]^\top$ and $\vec{\hat{q}}^{\vec{x}}_n$ is a realization of $\hat{Q}^{\vec{x}}_n$.} Then using \eqref{Eq: ACEP-UB}, $\ES{\norm{\vec{\epsilon}_k^{\vec{\lambda}}}{2}^2}$ can be lower bounded as
\begin{eqnarray}
\ES{\norm{\vec{\epsilon}_k^{\vec{\lambda}}}{2}^2}&\geq &\frac{\e{1-\frac{1}{N}}}{{2\pi \e{}}}\e{\frac{2}{N}\ES{\CENT{\vec{\epsilon}_k^{\vec{\lambda}}}{\mathcal{D}^{\vec{x}}_{k-1}}}}\nonumber\\
&=&\frac{\e{1-\frac{1}{N}}}{{2\pi \e{}}}\e{\frac{2}{N}\ES{\CENT{\vec{\lambda}_0}{\mathcal{D}^{\vec{x}}_{k-1}}}},
\end{eqnarray}
where the equality follows from translation invariance of the differential entropy and the update equation for the dual variables is linear in $\vec{\lambda}_{n}$. Using Lemma \ref{Lem: NE} from Appendix \ref{App: DDE}, $\ES{\CENT{\vec{\lambda}_{0}}{\mathcal{D}^{\vec{x}}_{k-1}}}$ can be lower bounded as
\begin{align}
&\ES{\CENT{\vec{\lambda}_{0}}{\mathcal{D}^{\vec{x}}_{k-1}}}\geq \mathsf{h}\left[\vec{\lambda}_{0}\right]-\sum_{t=0}^{k-1}\paren{\sum_{i=1}^M\log\abs{\mathcal{A}^{\vec{x}}_{i,t}}},
\end{align}
Thus, we have
\begin{align}
&\liminf_{k\longrightarrow\infty}\frac{1}{k}\log\ES{\norm{\vec{\epsilon}^{\vec{\lambda}}_k}{2}^2}\geq -\frac{2}{N} R_{\vec{x}}.
\end{align}
%%%%%%%%%%%%%%%%%%%%%%%%%%%%%%%%%%%%%%%%%%%%%%%%%%%%%%%%%%%%%%%%%%%%%%%%%%%%%%%%%%%%%%%%%%%%%%%%%%%%%%%%%%%%%%%%%%%%%%%%%%%%%%%%%%

\section{}\label{App: EDE-New}
\textcolor{black}{In this appendix, first, we establish a series of preliminary results in Subsection \ref{Sub-sec: Prelim}. Then, in Subsection \ref{Sub-sec: Proof}, we  use these preliminary results to prove Theorem \ref{Theo: EDE-New}. }
\subsection{Preliminary Lemmas}\label{Sub-sec: Prelim}
\textcolor{black}{ The next lemma shows that the quantization cells and their images under the quantized update rule are in the form of hypercubes. }
\begin{lemma}\label{Lem: Hyper-Cube}
\textcolor{black}{ The quantization cell $C_k\paren{\tilde{q}_k}$ is a hypercube in $\R^{N+M}$ for all $k$ and $\tilde{q}_k$. Also, \QM{C_k}{\tilde{q}_k}, \emph{i.e.,} the image of $C_k\paren{\tilde{q}_k}$ under the quantized update rule, is a hypercube.}
\end{lemma}
\begin{IEEEproof}
\textcolor{black}{ The proof of this lemma is based on mathematical induction. First, note that the quantized PD algorithm can be written as 
\begin{align}\label{EQ: QUR}
x^{i}_{k}&=T_{ii}x^{i}_{k-1}+\sum_{j=M+1}^{M+N}T_{ij}Q^{\lambda^j}_k+\mu c_i,\nonumber\\
\lambda^{j}_{k}&=\lambda^j_{k-1}+\sum_{i=1}^MT_{ji}Q^{x^i}_k-\mu b_j
\end{align}
 Since all the primal and dual variables are separately quantized, the quantization cells at time $k=0$, \emph{i.e.,} $C_0\paren{\tilde{q}_0}$s, are in the form of hypercubes. According to \eqref{EQ: QUR}, $x^{i}_{1}$ and $\lambda^{j}_{1}$ only depend on $x^{i}_{0}$ and $\lambda^{j}_{0}$, respectively, given $C_0\paren{\tilde{q}_0}$ (as $Q^{\vec{x}}_0$ and $Q^{\vec{\lambda}}_0$ are fixed given $C_0\paren{\tilde{q}_0}$). This observation implies that the image of $C_0\paren{\tilde{q}_0}$ under $\QM{\cdot}{\tilde{q}_0}$ is a hypercube. Now, assume that at time $k$ $C_k\paren{\tilde{q}_k}$ and $\QM{C_k}{\tilde{q}_k}$ are hypercubes. Since the quantization scheme is zoom-in, $\QM{C_k}{\tilde{q}_k}$ is quantized at time $k+1$. Therefore, $C_{k+1}\paren{\tilde{q}_{k+1}}$ is a hypercube as the primal/dual variables are independently quantized. Also, $Q^{\vec{x}}_{k+1}$ and $Q^{\vec{\lambda}}_{k+1}$ are fixed given $C_{k+1}\paren{\tilde{q}_{k+1}}$. Thus, the image of $C_{k+1}\paren{\tilde{q}_{k+1}}$ under $\QM{\cdot}{\tilde{q}_{k+1}}$ is a hypercube which completes the proof. }
\end{IEEEproof}
\textcolor{black}{The next lemma establishes an upper bound on the probability density function (PDF) of $\vec{y}_k$, \emph{i.e.,} $p_{\vec{y}_{k+1}}\paren{\vec{y}}$.}
\begin{lemma}\label{Lem: DUB}
\textcolor{black}{Consider the quantization cell $C_k\paren{\tilde{q}_k}$ at time $k$ where $\tilde{q}_k$ is a realization of $\tilde{Q}_k$. Let $G\paren{{\tilde{q}_{k}}}$ be the number of  quantization cells at time $k$ which their images under the quantized update rule overlap with that of $C_k\paren{\tilde{q}_k}$. Let $G^\star_{k+1}=\max_{\tilde{q}_k}G\paren{{\tilde{q}_k}}$. Then, the PDF of $\vec{y}_{k+1}$ can be upper bounded as 
\begin{align}
&p_{\vec{y}_{k+1}}\paren{\vec{y}}\leq \frac{P_{\rm max}G_{k+1}^\star}{\paren{\prod_{i=1}^{M+N}\abs{T_{ii}}}^{k+1}}
\end{align}
where $T_{ii}$ is the $i$th diagonal entry of $T$ and $P_{\rm max}$ is the maximum of PDF of $\vec{y}_0$.}
\end{lemma}
\begin{IEEEproof} 
\textcolor{black}{Let ${\rm B}\paren{\vec{y},r}$ be the hypercube centered at $\vec{y}$ with side length $r$.  Using Lebesgue differentiation theorem, the PDF of $\vec{y}_{k+1}$ can be written as
\begin{eqnarray}
p_{\vec{y}_{k+1}}\paren{\vec{y}}=\lim_{r\downarrow0}\frac{\PR{\vec{y}_{k+1}\in{\rm B}\paren{\vec{y},r}}}{{\rm Vol}\paren{{\rm B}\paren{\vec{y},r}}}
\end{eqnarray}
where ${\rm Vol}\paren{{\rm B}\paren{\vec{y},r}}$ is the volume of ${\rm B}\paren{\vec{y},r}$. $\PR{\vec{y}_{k+1}\in{\rm B}\paren{\vec{y},r}}$ can be written as 
\begin{align}\label{Eq: Prob-Expansion}
\PR{\vec{y}_{k+1}\in{\rm B}\paren{\vec{y},r}}=\sum_{\tilde{q}_k}\PR{\vec{y}_{k+1}\in{\rm B}\paren{\vec{y},r},\tilde{Q}_k=\tilde{q}_k}
\end{align}
where $\tilde{q}_k$ is a realization of $\tilde{Q}_k$. With a slight abuse of notation, let $\QM{C_k}{\tilde{q}_k}$ denote the image of $C_k\paren{\tilde{q}_k}$ under $\QM{\cdot}{\tilde{q}_k}$. Then, $\PR{\vec{y}_{k+1}\in{\rm B}\paren{\vec{y},r},\tilde{Q}_k=\tilde{q}_k}$ can be written as 
\begin{align}\label{Eq: Prob-Eq}
\PR{\vec{y}_{k+1}\in{\rm B}\paren{\vec{y},r},\tilde{Q}_k=\tilde{q}_k}&\nonumber\\
&\hspace{-5cm}=\PR{\vec{y}_{k+1}\in\paren{{\rm B}\paren{\vec{y},r}\cap \QM{C_k}{\tilde{q}_k}}}\IC{{\rm B}\paren{\vec{y},r}\cap \QM{C_k}{\tilde{q}_k}}
\end{align}
where the constant $\IC{{\rm B}\paren{\vec{y},r}\cap \QM{C_k}{\tilde{q}_k}}$ is equal to one if the intersection of $ {\rm B}\paren{\vec{y},r}$ and $\QM{C_k}{\tilde{q}_k}$ is non-empty and is equal to zero otherwise. % defined as 
%\begin{align}
%\IC{{\rm B}\paren{\vec{y},r}\cap \QM{C_k}{\tilde{q}_k}}=
%\left\{
%\begin{array}{cc}
%1 &  {\rm B}\paren{\vec{y},r}\cap \QM{C_k}{\tilde{q}_k}\neq\emptyset \nonumber\\
%0& {\rm B}\paren{\vec{y},r}\cap \QM{C_k}{\tilde{q}_k}=\emptyset\nonumber 
%\end{array}
%\right.\nonumber
%\end{align}
}

\textcolor{black}{Now, we find the set of initial conditions which allow the quantized PD algorithm to arrive in ${\rm B}\paren{\vec{y},r}\cap \QM{C_k}{\tilde{q}_k}$ at time $k+1$ with $\tilde{Q}_k=\tilde{q}_k$. Let $S_{k+1}={\rm B}\paren{\vec{y},r}\cap \QM{C_k}{\tilde{q}_k}$. Since $S_{k+1}$ is a subset of $\QM{C_k}{\tilde{q}_k}$, we can find $S_k\subset C_k\paren{\tilde{q}_k}$ such that $S_{k+1}=\QM{S_k}{\tilde{q}_k}$. Since the quantization scheme is zoom-in, $ C_k\paren{\tilde{q}_k}\subset \QM{C_{k-1}}{\tilde{q}_{k-1}}$ which implies $S_k\subset\QM{C_{k-1}}{\tilde{q}_{k-1}}$. Hence, we can find $S_{k-1}\subset C_{k-1}\paren{\tilde{q}_{k-1}}$ such that $S_k=\QM{S_{k-1}}{\tilde{q}_{k-1}}$. Using  backward induction, we can find a sequence of sets $\left\{S_n\right\}_{n=0}^k$ with $S_n\subset C_n\paren{\tilde{q}_n}$ such that we have $S_{n+1}=\QM{S_n}{\tilde{q}_n}$ with $\tilde{Q}_k=\tilde{q}_k$. %  Note that the $\QM{\cdot}{\tilde{q}_n}$ is one-to-one on $C_n\paren{\tilde{q}_n}$. Therefore, for any $\vec{y}_{k+1}\in \paren{{\rm B}\paren{\vec{y},r}\cap \QM{C_k}{\tilde{q}_k}}$, we can find $\vec{y}_0\in S_0$ such that if the quantized PD starts at $\vec{y}_0$, the output of the algorithm at time $k+1$ will be $\vec{y}_{k+1}$ and $\tilde{Q}_k=\tilde{q}_k$. 
Thus, we have 
\begin{align}\label{Eq: Prob-UpB}
\PR{\vec{y}_{k+1}\in\paren{{\rm B}\paren{\vec{y},r}\cap \QM{C_k}{\tilde{q}_k}}}&=\PR{\vec{y}_{0}\in S_0}\nonumber\\
%&\PR{\vec{y}_{0}\in S_0, \tilde{q}_k}\nonumber\\
&\leq P_{\rm max} {\rm Vol}\paren{S_0} 
\end{align}
where $P_{\rm max}$ is the maximum of PDF of $\vec{y}_0$ and ${\rm Vol}\paren{S_0}$ is the volume of $S_0$.}

\textcolor{black}{According to Lemma \ref{Lem: Hyper-Cube}, $\QM{C_k}{\tilde{q}_k}$ is a hypercube. Also, ${\rm B}\paren{\vec{y},r}$ is a hypercube, thus,  $S_{k+1}={\rm B}\paren{\vec{y},r}\cap \QM{C_k}{\tilde{q}_k}$ is a hypercube. Since $S_k\subset C_k\paren{\tilde{q}_k}$ and  $Q^{\vec{x}}_k$ and $Q^{\vec{\lambda}}_k$ are fixed given $C_k\paren{\tilde{q}_k}$,  $x^{i}_{k+1}$ ($\lambda^{j}_{k+1}$) only depends on $x^{i}_{k}$ ($\lambda^{j}_{k}$). Thus, $S_k$ is a hypercube. Also, the length of the $i$th side of $S_{k+1}$ is equal to the length of the $i$th side of $S_{k}$ multiplied by $\abs{T_{ii}}$. This is due to the fact that the quantized PD algorithm can be written as $x^{i}_{k}=T_{ii}x^{i}_{k-1}+\sum_{j=M+1}^{M+N}T_{ij}Q^{\lambda^j}_k+\mu c_i$ and $\lambda^{j}_{k}=\lambda^j_{k-1}+\sum_{i=1}^MT_{ji}Q^{x^i}_k-\mu b_j$. Hence, we have ${\rm Vol}\paren{S_{k+1}}={\rm Vol}\paren{S_k}\prod_{i=1}^{N+M}\abs{T_{ii}}$. Using backward induction, we have 
\begin{align}\label{Eq: Vol-Eq}
{\rm Vol}\paren{S_{k+1}}=\paren{\prod_{i=1}^{N+M}\abs{T_{ii}}}^{k+1}{\rm Vol}\paren{S_0}
\end{align}
The volume of $S_{k+1}$ can be upper bounded as
\begin{align}\label{Eq: Vol-Ineq}
{\rm Vol}\paren{S_{k+1}}&={\rm Vol}\paren{{\rm B}\paren{\vec{y},r}\cap \QM{C_k}{\tilde{q}_k}}\nonumber\\
&\leq {\rm Vol}\paren{{\rm B}\paren{\vec{y},r}}.
\end{align}
Combining \eqref{Eq: Prob-Expansion}, \eqref{Eq: Prob-Eq}, \eqref{Eq: Prob-UpB}, \eqref{Eq: Vol-Eq} and \eqref{Eq: Vol-Ineq}, we have
\begin{align}
&\PR{\vec{y}_{k+1}\in{\rm B}\paren{\vec{y},r}}\nonumber\\
&\leq \frac{P_{\rm max}}{\paren{\prod_{i=1}^{N+M}\abs{T_{ii}}}^{k+1}}\sum_{\tilde{q}_k}\IC{{\rm B}\paren{\vec{y},r}\cap \QM{C_k}{\tilde{q}_k}}
\end{align}
Note that $\lim_{r\downarrow0}\IC{{\rm B}\paren{\vec{y},r}\cap \QM{C_k}{\tilde{q}_k}}=\IC{\vec{y}\in \QM{C_k}{\tilde{q}_k}}$ where $\IC{\vec{y}\in \QM{C_k}{\tilde{q}_k}}$ is equal to one  if $\vec{y}\in \QM{C_k}{\tilde{q}_k}$ and is equal to zero otherwise. Let $G_{k+1}^{\vec{y}}=\sum_{\tilde{q}_k}\IC{\vec{y}\in \QM{C_k}{\tilde{q}_k}}$, \emph{i.e.,} the number of quantization cells at time $k$ which their images under the quantized update rule \emph{jointly} overlap at $\vec{y}$. Thus, we have 
\begin{align}\label{}
&p_{\vec{y}_{k+1}}\paren{\vec{y}}\leq \frac{P_{\rm max}G_{k+1}^{\vec{y}}}{\paren{\prod_{i=1}^{N+M}\abs{T_{ii}}}^{k+1}}
\end{align} }

\textcolor{black}{ Note that $\max_{\vec{y}}G_{k+1}^{\vec{y}}$ is the maximum number of quantization cells at time $k$ which their images jointly overlap with each other under $\QM{\cdot}{\tilde{q}_k}$. For a quantization cell $C_k\paren{\tilde{q}_k}$, let $G^\prime\paren{{\tilde{q}_k}}$ be the maximum number of quantization cells at time $k$ which their images under  the quantized update rule \emph{jointly} overlap with the image of  $C_k\paren{\tilde{q}_k}$. Also, let $G\paren{{\tilde{q}_k}}$ be the number of  quantization cells at time $k$ which their images under  the quantized update rule overlap with the image of  $C_k\paren{\tilde{q}_k}$. Clearly, we have $G^\prime\paren{{\tilde{q}_k}}\leq G\paren{{\tilde{q}_k}}$. Thus, we have 
\begin{align}
\max_{\vec{y}}G^{k+1}_{\vec{y}}&= \max_{{\tilde{q}_k}}G^\prime\paren{{\tilde{q}_k}}\nonumber\\
&\leq\max_{{\tilde{q}_k}}G\paren{{\tilde{q}_k}}\nonumber\\
&=G^\star_{k+1}
\end{align}
Hence, the PDF of $\vec{y}_{k+1}$ can be upper bounded as 
\begin{align}
&p_{\vec{y}_{k+1}}\paren{\vec{y}}\leq \frac{P_{\rm max}G_{k+1}^\star}{\paren{\prod_{i=1}^{N+M}\abs{T_{ii}}}^{k+1}}
\end{align} } 
\end{IEEEproof} 
\textcolor{black}{The next two results are used in Lemma \ref{Lem: NLP} to obtain an upper bound on $G^\star_{k+1}$. The next lemma derives a necessary condition for two hypercubes to overlap.}
\begin{lemma}\label{Lem: HQC}
\textcolor{black}{Consider two hypercubes $C_1$ and $C_2$ in $\R^{N+M}$ which the length of the $j$th side of $C_i$ is given by $l_i^j>0$ for $i=1,2$ and $1\leq j\leq N+M$. Let $\vec{y}_{C_1}$ and $\vec{y}_{C_2}$ be any two points in $C_1$ and $C_2$, respectively. Then, a necessary condition for $C_1$ and $C_2$  to overlap with each other is given by 
\begin{align}
-\paren{l^j_1+l^j_2}\leq \vec{y}^j_{C_2}-\vec{y}^j_{C_1}\leq l^j_1+l^j_2 \quad \textit{for all $j$}\nonumber
\end{align}
for any $\vec{y}_{C_1}\in C_1$ and $\vec{y}_{C_2}\in C_2$ where $\vec{y}^j_{C_i}$ is the $j$th entry of  $\vec{y}_{C_i}$.}
 \end{lemma}
\begin{IEEEproof}
\textcolor{black}{Without loss of generality, we assume  $C_1$ and $C_2$ which are specified by $\vec{y}^{\rm min}_{C_1}\leq \vec{y}\leq \vec{y}^{\rm max}_{C_1}$ and $\vec{y}^{\rm min}_{C_2}\leq \vec{y}\leq \vec{y}^{\rm max}_{C_2}$, receptively, where $\vec{y}\in\R^{N+M}$. Then, for $\vec{y}\in C_1\cap C_2$, we have 
\begin{align}
\vec{y}^{\rm min}_{C_1}-\vec{y}_{C_1}&\leq \vec{y}-\vec{y}_{C_1}\leq \vec{y}^{\rm max}_{C_1}-\vec{y}_{C_1}\nonumber\\
\vec{y}^{\rm min}_{C_2}-\vec{y}_{C_2}&\leq \vec{y}-\vec{y}_{C_2}\leq \vec{y}^{\rm max}_{C_2}-\vec{y}_{C_2}\nonumber
\end{align}
where $\leq$ indicates component-wise inequality. Combining the above inequalities, we obtain a necessary condition for $C_1$ and $C_2$ to overlap as follows 
\begin{align}
-\paren{\vec{y}^{\rm max}_{C_2}-\vec{y}_{C_2}}-&\paren{\vec{y}_{C_1}-\vec{y}^{\rm min}_{C_1}}\nonumber\\
&\leq \vec{y}_{C_2}-\vec{y}_{C_1}\leq \nonumber\\
&\paren{\vec{y}_{C_2}-\vec{y}^{\rm min}_{C_2}}+\paren{\vec{y}^{\rm max}_{C_1}-\vec{y}_{C_1}}\nonumber
\end{align}
 Note that $\vec{y}^{\rm max}_{C_i}-\vec{y}_{C_i}\leq \vec{y}^{\rm max}_{C_i}-\vec{y}^{\rm min}_{C_i}$ and $\vec{y}_{C_i}-\vec{y}^{\rm min}_{C_i}\leq \vec{y}^{\rm max}_{C_i}-\vec{y}^{\rm min}_{C_i}$. Thus, the above condition can be further relaxed to the following necessary condition for $C_1$ and $C_2$ to overlap:
\begin{align}
-\paren{\vec{y}^{\rm max}_{C_2}-\vec{y}^{\rm min}_{C_2}}-&\paren{\vec{y}^{\rm max}_{C_1}-\vec{y}^{\rm min}_{C_1}}\nonumber\\
&\leq \vec{y}_{C_2}-\vec{y}_{C_1}\leq \nonumber\\
&\paren{\vec{y}^{\rm max}_{C_2}-\vec{y}^{\rm min}_{C_2}}+\paren{\vec{y}^{\rm max}_{C_1}-\vec{y}^{\rm min}_{C_1}}
\end{align}
The proof is complete by appealing to the fact $l^j_i$ is the $j$th element of the $\vec{y}^{\rm max}_{C_i}-\vec{y}^{\rm min}_{C_i}$.}
\end{IEEEproof}
\textcolor{black}{Consider the lattice $\Lambda\paren{\vec{y}}=\vec{y}+\delta^{\rm min}_k\Z^{N+M}$ where, $\vec{y}\in \R^{N+M}$, $\delta^{\rm min}_k$ is the minimum quantization step at time $k$, and $\Z^{N+M}$ is the  $N+M$ dimensional integer lattice in $\R^{N+M}$. Next, we define the minimum distance assignment rule which assigns a unique lattice point of $\Lambda\paren{\vec{y}}$ to each quantization cell. Then, in Lemma \ref{Lem: NLP-aux}, we derive an upper bound on the number of quantization cells which a lattice point of $\Lambda\paren{\vec{y}}$ can be assigned to under minimum distance assignment rule. }
\begin{definition}
\textcolor{black}{Consider the lattice $\Lambda\paren{\vec{y}}$. Under minimum distance assignment rule, a lattice point of $\Lambda\paren{\vec{y}}$ is assigned to each quantization cell as follows. If a quantization cell contains only one point of $\delta^{\rm min}_k$, that point is assigned to the corresponding cell. If a quantization cell contains more than one lattice point, then, a lattice point with the smallest distance to its cell representative is assigned to the corresponding cell.}
\end{definition}

\textcolor{black}{Since each side length of quantization cells at time $k$ is greater or equal to $\delta^{\rm min}_k$, every quantization cell at time $k$ at least contains one point of $\Lambda\paren{\vec{y}}$. 
\begin{lemma}\label{Lem: NLP-aux}
Every lattice point in $\Lambda\paren{\vec{y}}$ can at most  be assigned to $G^\star_k$ quantization cells at time $k$.
\end{lemma}}
\begin{IEEEproof}
\textcolor{black}{This lemma is proved by contradiction. Let $G\paren{{\tilde{q}_{k-1}}}$ be the number of  quantization cells at time $k-1$ which their images under the quantized update rule overlap with that of $C_{k-1}\paren{\tilde{q}_{k-1}}$. Then, $G^\star_{k}$ can be written as $G^\star_{k}=\max_{\tilde{q}_{k-1}}G\paren{{\tilde{q}_{k-1}}}$. Assume that there exists a point $\vec{y}+\delta^{\rm min}_k\vec{I}$ with $\vec{I}\in\Z^{N+M}$ which can be assigned to $G>G^\star_k$ quantization cells at time $k$. Since, under a zoom-in quantization scheme, at time $k$, the image of a quantization cell at time $k-1$ is quantized, $\vec{y}+\delta^{\rm min}_k\vec{I}$ belongs to the images of $G$ quantization cells at time $k-1$. This observation implies that image of $G$ quantization cells at time $k-1$ under the quantized update rule overlap with each other at $\vec{y}+\delta^{\rm min}_k\vec{I}$. Thus, we have  $G\leq G^\star_k$ which contradicts with our assumption. }
\end{IEEEproof}
\begin{lemma}\label{Lem: NLP}
\textcolor{black}{Let $\beta_T$ be the number points in the lattice $T\Z^{N+M}$ which lie in a hypercube centered at the origin with the $i$th side length equal to $4\rho\abs{T_{ii}}+2\norm{T}{\infty}$ where $\norm{\cdot}{\infty}$ denotes the norm infinity, $T_{ii}$ is the $i$th diagonal entry of matrix $T$, and $\Z^{N+M}$ is the lattice of integers in $\R^{N+M}$. %Let $\beta_T$ be the number of lattice points $\vec{I}=\left[I_1,\cdots,I_{N+M}\right]^\top\in\Z^{M+n}$ which satisfy the following system of inequalities 
%\begin{align}
%-2\frac{\delta^{\rm max}_k}{\delta^{\rm min}_k}\abs{T_{jj}}-\norm{T}{\infty}\leq\sum_{m}T_{jm}{I}_m\leq 2\frac{\delta^{\rm max}_k}{\delta^{\rm min}_k}\abs{T_{jj}}+\norm{T}{\infty} \quad \forall j\nonumber
%\end{align}
Then, $G^\star_{k+1}\leq \beta_T^{k+1}$.}
\end{lemma}
\begin{IEEEproof}
\textcolor{black}{Consider two distinct quantization cells at time $k$ $C_k\paren{\tilde{q}_k}$ and $C_k\paren{\tilde{q}^\prime_k}$ with the cell representatives $\vec{y}_r\paren{\tilde{q}_k}\in C_k\paren{\tilde{q}_k}$ and $\vec{y}_r\paren{\tilde{q}^\prime_k}\in C_k\paren{\tilde{q}^\prime_k}$. Let $\QM{C_k\paren{\tilde{q}_k}}{\tilde{q}_k}$ and $\QM{C_k\paren{\tilde{q}_k^\prime}}{\tilde{q}_k^\prime}$ be the images of $C_k\paren{\tilde{q}_k}$ and $C_k\paren{\tilde{q}_k^\prime}$, respectively,  under the quantized update rule. $\QM{C_k\paren{\tilde{q}_k}}{\tilde{q}_k}$ and $\QM{C_k\paren{\tilde{q}_k^\prime}}{\tilde{q}_k^\prime}$ are hypercubes in $\R^{N+M}$, thus, Lemma \ref{Lem: HQC} can be used to obtain a necessary condition for $\QM{C_k\paren{\tilde{q}_k}}{\tilde{q}_k}$ and $\QM{C_k\paren{\tilde{q}_k^\prime}}{\tilde{q}_k^\prime}$ to overlap. Note that the $j$th side length of $\QM{C_k\paren{\tilde{q}_k}}{\tilde{q}_k}$ is less than or equal to $\abs{T_{jj}}\delta^{\rm max}_k$. This due to the facts that $j$th side length of $\QM{C_k\paren{\tilde{q}_k}}{\tilde{q}_k}$ is equal to the $j$th side length of $C_k\paren{\tilde{q}_k}$ multiplied by $\abs{T_{jj}}$ and each side length of quantization cells at time $k$ is less than or equal to $\delta^{\rm max}_k$. Similarly, the $j$th side length of $T\paren{C_k\paren{\tilde{q}^\prime_k},\tilde{q}^\prime_k}$ is less than or equal to $\abs{T_{jj}}\delta^{\rm max}_k$. %Substituting $C_1=\QM{C_k\paren{\tilde{q}_k}}{\tilde{q}_k}$, $C_2=\QM{C_k\paren{\tilde{q}^\prime_k}}{\tilde{q}^\prime_k}$, $\vec{y}_{C_1}=\QM{\vec{y}_r\paren{\tilde{q}_k}}{\tilde{q}_k}$, $\vec{y}_{C_2}=\QM{\vec{y}_r\paren{\tilde{q}^\prime_k}}{\tilde{q}^\prime_k}$, $l^j_i=\abs{T_{jj}}\delta^{\rm max}_k$ ($i=1,2$), in 
Thus, using Lemma \ref{Lem: HQC}, the necessary condition for $\QM{C_k\paren{\tilde{q}_k}}{\tilde{q}_k}$ and $\QM{C_k\paren{\tilde{q}^\prime_k}}{\tilde{q}^\prime_k}$ to overlap is given by 
\begin{align}\label{Eq: ini-nec-con}
-2\delta^{\rm max}_k\abs{T_{jj}}\leq\hat{T}_j\paren{\vec{y}_r\paren{\tilde{q}_k},\tilde{q}_k}&-\hat{T}_j\paren{\vec{y}_r\paren{\tilde{q}^\prime_k},\tilde{q}^\prime_k}\nonumber\\
& \leq 2\delta^{\rm max}_k\abs{T_{jj}}  \quad \forall j
\end{align}
Note that the quantization does not have any impact on the  representative of quantization cells. Thus, the quantized update rule for cell representatives is the same as the unquantized update rule. Thus, $\QM{\vec{y}_r\paren{\tilde{q}_k}}{\tilde{q}_k}-\QM{\vec{y}_r\paren{\tilde{q}^\prime_k}}{\tilde{q}^\prime_k}=T\paren{\vec{y}_r\paren{\tilde{q}_k}-\vec{y}_r\paren{\tilde{q}^\prime_k}}$.
%,   the output of the quantized PD update rule for a primal (dual) variable  at time $k+1$ is independent of the values of other variables (as the quantized versions of other variables are fixed on $C_k\paren{\tilde{q}_k}$).  Moreover, the length of the $i$th side of $\QM{C_k}{\tilde{q}_k}$ is equal to the length of the $i$th side of $C_k\paren{\tilde{q}_k}$ multiplied by $T_{ii}$. Thus, we have ${\rm Vol}\paren{C_k}=\prod_{i=1}^{N+M}\abs{T_{ii}}{\rm Vol}\paren{\QM{C_k}{\tilde{q}_k}}$ and 
% Thus, we have 
%\begin{align}
%-2\delta^{\rm max}_k\abs{T_{jj}}\leq \sum_i T_{ji}\paren{\vec{y}^i_r\paren{\tilde{q}_k}-\vec{y}^i_r\paren{\tilde{q}^\prime_k}} \leq 2\delta^{\rm max}_k\abs{T_{jj}}\nonumber
%\end{align}
}
\textcolor{black}{ The cell representative of $C_k\paren{\tilde{q}^\prime_k}$ can be written as $\vec{y}_r\paren{\tilde{q}^\prime_k}=\vec{y}_r\paren{\tilde{q}_k}+\delta^{\rm min}_k\vec{I}+\vec{\psi}$ where $\vec{y}_r\paren{\tilde{q}_k}+\delta^{\rm min}_k\vec{I}$ is the lattice point of $\Lambda\paren{\vec{y}_r\paren{\tilde{q}_k}}$  assigned to $C_k\paren{\tilde{q}^\prime_k}$ and $\vec{\psi}$ is an $N+M$ dimensional vector with the $j$th entry satisfying $-\delta^{\rm min}_k< \psi_j<\delta^{\rm min}_k$. Thus, we have $\QM{\vec{y}_r\paren{\tilde{q}_k}}{\tilde{q}_k}-\QM{\vec{y}_r\paren{\tilde{q}^\prime_k}}{\tilde{q}^\prime_k}=\delta^{\rm min}_kT\vec{I}+T\vec{\psi}$.
Using \eqref{Eq: ini-nec-con}, a necessary condition for $\QM{C_k\paren{\tilde{q}_k}}{\tilde{q}_k}$ and $\QM{C_k\paren{\tilde{q}^\prime_k}}{\tilde{q}^\prime_k}$ to overlap can be obtained as 
\begin{align}
-2\rho T_d-\frac{1}{\delta^{\rm min}_k}T\vec{\psi}\leq T\vec{I}\leq2\rho T_d-\frac{1}{\delta^{\rm min}_k}T\vec{\psi}\nonumber
\end{align}
where $T_d=\left[\abs{T_{ii}}\right]_i^\top$ and $\rho=\frac{\delta^{\rm max}_k}{\delta^{\rm min}_k}$ and $\leq$ indicates component-wise inequality. Using the fact that $\abs{\frac{{\psi}_m}{\delta^{\rm min}_k}}\leq 1$, the above condition can be relaxed to 
\begin{align}\label{Eq: Overlap-Cond}
-2\rho T_d-\norm{T}{\infty}\1_{N+M}\leq T\vec{I}\leq2\rho T_d+\norm{T}{\infty}\1_{N+M}
\end{align}
 where $\1_{N+M}$ is an $N+M$ dimensional vector with all entries equal to one. Recall that $G\paren{\tilde{q}_k}$ is the number of quantization cells which overlap with $C_k\paren{\tilde{q}_k}$. The number of quantization cells satisfying \eqref{Eq: Overlap-Cond}  provides an upper bound on $G\paren{\tilde{q}_k}$ as  \eqref{Eq: Overlap-Cond} is a necessary condition for two cells to overlap with each other. Let $\beta_T$ be the number of vectors $\vec{I}\in\Z^{N+M}$ which satisfy \eqref{Eq: Overlap-Cond}. According to Lemma \ref{Lem: NLP-aux},  for any $\vec{I}\in\Z^{N+M}$,  the lattice point $\vec{y}_r\paren{\tilde{q}_k}+\delta^{\rm min}_k\vec{I}$ can be assigned to at most $G^\star_k$  quantization cells at time $k$. Thus, the number of quantization cell for which \eqref{Eq: Overlap-Cond} holds is upper bounded by $\beta_T G^\star_k$ and, we have $G\paren{\tilde{q}_k}\leq \beta_T G^\star_k$. % For every vector $\vec{I}\in\Z^{N+M}$ satisfying \eqref{Eq: Overlap-Cond}, the quantization cells corresponding to $\vec{y}_r\paren{\tilde{q}_k}+\delta^{\rm min}_k\vec{I}$ may overlap with  $C_k\paren{\tilde{q}_k}$ as \eqref{Eq: Overlap-Cond} is a necessary condition for two cells to overlap with each other.  herefore, if a quantization cell at time $k$ overlaps with $C_k\paren{\tilde{q}_k}$, the vector $\vec{I}$, corresponding to its assigned lattice point, satisfies \eqref{Eq: Overlap-Cond} Since, . %Note that any cell representative at time $k$ can be written as $\vec{y}_r\paren{\tilde{q}_k}+\delta^{\rm min}_k\vec{I}+\vec{\psi}$ for some $\vec{I}\in\Z^{N+M}$ and $\vec{\psi}$.  Thus, every solution of the above system of inequalities corresponds to at most $G^\star_k$ quantization cells which may overlap with $C_k\paren{\tilde{q}_k}$. Then, the total number of quantization cells which overlap with $C_k\paren{\tilde{q}_k}$. 
		Since this bound is independent of $C_k\paren{\tilde{q}_k}$, we have $G^\star_{k+1}\leq \beta_T G^\star_k$. Following a similar argument, it can be easily shown than $G^\star_1\leq \beta_T$. Hence, we have $G^\star_{k+1}\leq \beta_T^{k+1}$ which completes the proof. 
}
\end{IEEEproof}
%%%%%%%%%%%%%%%%%%%%%%%%%%%%%%%%%%%%%%%%%%%%%%%%%%%%%%%%%%%%%%%%%%%%%%%%%%%%%%%%%%%%%%%%%%%%%%%%%%%%%%%%%%%%%%%%%%%%
\subsection{Proof of Theorem \ref{Theo: EDE-New}}\label{Sub-sec: Proof}
\textcolor{black}{Now, we are ready to prove Theorem \ref{Theo: EDE-New}. To this end, first, we use Lemma \ref{Lem: DUB}  to establish a lower bound on the differential entropy of $\vec{y}_{k+1}$ a follows:
\begin{align}\label{Eq: Ent-LB}
\ENT{\vec{y}_{k+1}}%&=-\int \logp{ p_{\vec{y}_{k+1}}\paren{\vec{y}}}p_{\vec{y}_{k+1}}\paren{\vec{y}}d\vec{y}\nonumber\\
&=\int -\logp{ p_{\vec{y}_{k+1}}\paren{\vec{y}}}p_{\vec{y}_{k+1}}\paren{\vec{y}}d\vec{y}\nonumber\\
&\stackrel{(a)}{\geq} -\log{ \frac{P_{\rm max}G_{k+1}^\star}{\paren{\prod_{i=1}^{N+M}\abs{T_{ii}}}^{k+1}}} \int p_{\vec{y}_{k+1}}\paren{\vec{y}}d\vec{y}\nonumber\\
&= -\log{ \frac{P_{\rm max}G_{k+1}^\star}{\paren{\prod_{i=1}^{N+M}\abs{T_{ii}}}^{k+1}}} 
\end{align}
where $(a)$ follows from the fact that $-\logp{x}$ is decreasing in $x$. Let $\epsilon_{k+1}=\vec{y}_{k+1}-\vec{y}^\star$. Substituting $A=\Omega$, where $\Omega$ is the sample space of the underlying probability space, and $\vec{z}=\epsilon_{k+1}$ in equation \eqref{Eq: CEP-UB}, we have  
\begin{align}\label{Eq: MSE-LB}
\ES{\norm{\vec{\epsilon}_{k+1}}{2}^2}%&\geq \frac{\e{1-\frac{1}{N+M}}}{{2\pi \e{}}}\e{\frac{2}{N+M}\ENT{\vec{\epsilon}_{k+1}}}\nonumber\\
&\geq  \frac{\e{1-\frac{1}{N+M}}}{{2\pi \e{}}}\e{\frac{2}{N+M}\ENT{\vec{y}_{k+1}}}
\end{align}
Combining \eqref{Eq: Ent-LB} and \eqref{Eq: MSE-LB}, we have 
\begin{align}\label{Eq: EXP-LB-A3}
\log{\ES{\norm{\vec{\epsilon}_{k+1}}{2}^2}}&\geq \log{\frac{\e{1-\frac{1}{N+M}}}{{2\pi \e{}}}}\nonumber\\
&-\frac{2}{N+M}\log{ \frac{P_{\rm max}G_{k+1}^\star}{\paren{\prod_{i=1}^{N+M}\abs{T_{ii}}}^{k+1}}} 
\end{align}}

\textcolor{black}{Combining \eqref{Eq: EXP-LB-A3} and Lemma \ref{Lem: NLP}, we have 
\begin{align}
\liminf_{k\rightarrow\infty}\frac{1}{k+1}&\log{\ES{\norm{\vec{\epsilon}_{k+1}}{2}^2}}\nonumber\\
&\geq -\frac{2}{N+M}\log{\frac{ \beta_T}{\paren{\prod_{i=1}^{N+M}\abs{T_{ii}}}}} 
\end{align}
which completes the proof.
}
%%%%%%%%%%%%%%%%%%%%%%%%%%%%%%%%%%%%%%%%%%%%%%%%%%%%%%%%%%%%%%%%%%%%%%%%%%%%%%%%%%%%%%%%%%%%%%%%%%%%%%%%%%%%%%%%%%%%%%%%%%%%%%%%%%%
\section{Proof of Theorem \ref{Theo: Conv}}\label{App:Conv}
\textcolor{black}{In this appendix, first, we show that $x^i_{k+1}$ and $\lambda^j_{k+1}$ belong to $I^{x^i}_{k+1}$ and $I^{\lambda^j}_{k+1}$, respectively. Note that $\abs{x^i_{k+1}-x^i_k-\left\lfloor \frac{x^i_{k+1}-x^i_k}{\delta_k}\right\rfloor\delta_k}\leq \delta_k$. Thus, 
\begin{align}
&\abs{x^i_{k+1}-C^{x^i}_{k+1}}&\nonumber\\
&= \abs{x^i_{k+1}-x^i_{k}+x^i_{k}-{\rm Q}_{{\rm a},k}\paren{x^i_k}-\left\lfloor \frac{x^i_{k+1}-x^i_k}{\delta_k}\right\rfloor\delta_k}\nonumber\\
&{\leq} \abs{x^i_{k+1}-x^i_{k}-\left\lfloor \frac{x^i_{k+1}-x^i_k}{\delta_k}\right\rfloor\delta_k}+\abs{x^i_{k}-{\rm Q}_{{\rm a},k}\paren{x^i_k}}\nonumber\\
&\stackrel{\paren{a}}{\leq2}\delta_k\nonumber\\
&\leq \left\lceil \frac{2}{\alpha}\right\rceil\alpha\delta_k\nonumber\\
&= \left\lceil \frac{2}{\alpha}\right\rceil\delta_{k+1}
\end{align}
where $(a)$ follows from the fact that the quantization error is less than the quantization step $\delta_k$. Thus, $x^i_{k+1}$ belongs to $I^{x^i}_{k+1}$. Similarly, we have $\abs{\lambda^j_{k+1}-\lambda^j_k-\left\lfloor \frac{\lambda^j_{k+1}-\lambda^j_{k}}{\delta_k}\right\rfloor\delta_k}\leq \delta_k$. Thus, we have 
\begin{align}
&\abs{\lambda^j_{k+1}-C^{\lambda^j}_{k+1}}\nonumber\\
&\leq \abs{\lambda^j_{k+1}-\lambda^j_k-\left\lfloor \frac{\lambda^j_{k+1}-\lambda^j_{k}}{\delta_k}\right\rfloor\delta_k}+\abs{\lambda^j_k-{\rm Q}_{{\rm a},k}\paren{\lambda^j_k}}\nonumber\\
&\leq 2\delta_k\nonumber\\
&\leq \left\lceil \frac{2}{\alpha}\right\rceil\delta_{k+1}
\end{align}
which implies that $\lambda^j_{k+1}$ belongs to $I^{\lambda^j}_{k+1}$.}

\textcolor{black}{Now, we prove the convergence of the PD algorithm under the quantization scheme $\mathcal{Q}_{\rm a}$. Note that the PD update rule can be written as 
\begin{align}
{x}^{i}_{k}&=x^{i}_{k-1}+\mu_{k-1} \paren{\frac{d }{d x^i}U_i\paren{x^i_{k-1}}-A^\top_i\vec{\lambda}_{k-1}},\nonumber  \\
&\hspace{2cm}+\mu_{k-1}A^\top_i\paren{\vec{\lambda}_{k-1}-Q^{\vec{\lambda}}_{k-1}}\nonumber\\
{\lambda}^{j}_{k}&=\lambda^j_{k-1}+\mu_{k-1}\paren{A_j\vec{x}_{k-1}-b_j}+\mu_{k-1}A_j\paren{Q^{\vec{x}}_{k-1}-\vec{x}_{k-1}}\nonumber
\end{align}
Since the unquantized update rule forms a contraction map, we have 
\begin{align}\label{Eq: Norm-UB}
&\norm{
\vec{\epsilon}_k}{}\leq\alpha \norm{\vec{\epsilon}_{k-1}}{}+ \mu_{k-1}
\norm{
\left[
\begin{array}{cc}
\vec{0}&A^\top\\
-A&\vec{0}
\end{array}
\right]\left[
\begin{array}{c}
\vec{x}_{k-1}-\quant{\vec{x}_{k-1}}{\vec{x}}{k-1}\\
\vec{\lambda}_{k-1}-\quant{\vec{\lambda}_{k-1}}{\vec{\lambda}}{k-1}
\end{array}
\right]}{}
\end{align}
where $\norm{\cdot}{}$ is the norm in which the unquantized PD algorithm forms a contraction mapping. The second term in the right hand side of \eqref{Eq: Norm-UB} can be upper bounded as 
\begin{align}\label{Eq: Norm-UB-1}
& \mu_{k-1}
\norm{
\left[
\begin{array}{cc}
\vec{0}&A^\top\\
-A&\vec{0}
\end{array}
\right]\left[
\begin{array}{c}
\vec{x}_{k-1}-\quant{\vec{x}_{k-1}}{\vec{x}}{k-1}\\
\vec{\lambda}_{k-1}-\quant{\vec{\lambda}_{k-1}}{\vec{\lambda}}{k-1}
\end{array}
\right]}{},\nonumber\\
&\stackrel{(a)}{\leq} \mu_{k-1}C
\norm{
\left[
\begin{array}{cc}
\vec{0}&A^\top\\
-A&\vec{0}
\end{array}
\right]\left[
\begin{array}{c}
\vec{x}_{k-1}-\quant{\vec{x}_{k-1}}{\vec{x}}{k-1}\\
\vec{\lambda}_{k-1}-\quant{\vec{\lambda}_{k-1}}{\vec{\lambda}}{k-1}
\end{array}
\right]}{\infty},\nonumber\\
&\stackrel{(b)}{\leq} \min_i\frac{1}{\abs{U^{\rm min}_i}}C
\norm{
\left[
\begin{array}{cc}
\vec{0}&A^\top\\
-A&\vec{0}
\end{array}
\right]}{\infty}\norm{\left[
\begin{array}{c}
\vec{x}_{k-1}-\quant{\vec{x}_{k-1}}{\vec{x}}{k-1}\\
\vec{\lambda}_{k-1}-\quant{\vec{\lambda}_{k-1}}{\vec{\lambda}}{k-1}
\end{array}
\right]}{\infty},\nonumber\\
&\stackrel{(c)}{=} \min_i\frac{1}{\abs{U^{\rm min}_i}}C
\norm{
\left[
\begin{array}{cc}
\vec{0}&A^\top\\
-A&\vec{0}
\end{array}
\right]}{\infty}\delta_{k-1},
\end{align}
where $C$ is a positive constant, $(a)$ follows from the fact that all the norms on a finite dimensional vector space are equivalent, $(b)$ follows from $0<\mu_{n}< \min_i\frac{1}{\abs{U^{\rm min}_i}}$, and the multiplicative property of norm infinity, and $(c)$ follows from the fact that the quantization error at time $k-1$ is less than or equal to $\delta_{k-1}$ for all primal and dual variables.}

\textcolor{black}{Using \eqref{Eq: Norm-UB} and \eqref{Eq: Norm-UB-1}, $\norm{\vec{\epsilon}_k}{}$ can be upper bounded as 
\begin{align}\label{Eq: Exp-Con}
&\norm{{\vec{\epsilon}_k}}{}\!\leq\! 
\alpha^k\!\norm{\vec{\epsilon}_{0}}{}\!+\!\min_i\frac{1}{\abs{U^{\rm min}_i}}C
\norm{
\left[
\begin{array}{cc}
\vec{0}&A^\top\\
-A&\vec{0}
\end{array}
\right]}{\infty}\!\sum_{i=1}^{k}\!\alpha^{i-1}\delta_{k-i}\nonumber\\
&\stackrel{(a)}{=} \alpha^k\norm{\vec{\epsilon}_{0}}{}+\min_i\frac{1}{\abs{U^{\rm min}_i}}C
\norm{
\left[
\begin{array}{cc}
\vec{0}&A^\top\\
-A&\vec{0}
\end{array}
\right]}{\infty}k\alpha^k
\end{align}
where $(a)$ follows from the fact that $\delta_k=\alpha^{k+1}$. The last inequality in \eqref{Eq: Exp-Con} implies that the PD algorithm under the quantization scheme $\mathcal{Q}_{\rm a }$ converges exponentially to the optimal solution.}
\bibliographystyle{IEEEtran}

\begin{thebibliography}{1}
\providecommand{\url}[1]{#1}
\csname url@samestyle\endcsname
\providecommand{\newblock}{\relax}
\providecommand{\bibinfo}[2]{#2}
\providecommand{\BIBentrySTDinterwordspacing}{\spaceskip=0pt\relax}
\providecommand{\BIBentryALTinterwordstretchfactor}{4}
\providecommand{\BIBentryALTinterwordspacing}{\spaceskip=\fontdimen2\font plus
\BIBentryALTinterwordstretchfactor\fontdimen3\font minus
  \fontdimen4\font\relax}
\providecommand{\BIBforeignlanguage}[2]{{%
\expandafter\ifx\csname l@#1\endcsname\relax
\typeout{** WARNING: IEEEtran.bst: No hyphenation pattern has been}%
\typeout{** loaded for the language `#1'. Using the pattern for}%
\typeout{** the default language instead.}%
\else
\language=\csname l@#1\endcsname
\fi
#2}}
\providecommand{\BIBdecl}{\relax}

\BIBdecl
\bibitem{KMT98}
F.~P.~Kelly, A.~Maulloo, and D.~Tan, ``Rate control for communication
networks: Shadow prices, proportional fairness and stability," \emph{J. Oper. Res}., vol. 49, no. 3, pp. 237–252,  1998.

\bibitem{SS07}
S.~Shakkottai and R.~Srikant. ``Network optimization and control.", \emph{Found. Trends Netw.}, vol.~2, no.~3 pp.~271-379, 2007.
%\bibitem{NO09}
%A.~Nedi\'{c} and A.~Ozdaglar, ``Distributed Subgradient Methods for Multi-Agent Optimization," \emph{IEEE Trans. Autom. Control}, vol.~54, no.~1, pp.~48-61, Jan. 2009.
\bibitem{AB05}
T.~Alpcan and T.~Basar, ``A Utility-Based Congestion Control Scheme for Internet-Style Networks with Delay," \emph{IEEE Trans. Netw.}, vol.~13, no.~6, pp.~1261-1274,  2005.

%\bibitem{NB01}
%A.~Nedi\'{c} and D.P.~Bertsekas, ``Incremental Subgradient Methods for Nondifferentiable Optimization", \emph{SIAM J. Control Optim.}, vol.~12, no.~1, pp.~109–138 2001.

%\bibitem{NO09}
%A.~Nedi\'{c} and A.~Ozdaglar, ``Distributed Subgradient Methods for Multi-Agent Optimization," \emph{IEEE Trans. Autom. Control}, vol.~54, no.~1, pp.~48-61, Jan. 2009.

%\bibitem{NOP10}
%A.~Nedi\'{c}, A.~Ozdaglar and P.A.~Parrilo, ``Constrained Consensus and Optimization in Multi-Agent Networks," \emph{IEEE Trans. Autom. Control}, vol.~55, no.~4, pp.~922-938, April 2010.


\bibitem{Nedic08}
A.~Nedi\'{c}, A.~Olshevsky, A.~Ozdaglar and J.N.~Tsitsiklis, ``Distributed subgradient methods and quantization effects," IEEE Conf. on Decision and Control Conference (CDC), pp.4177--4184,  2008.

\bibitem{Rabbat05}
M.G.~Rabbat and R.D.~Nowak, ``Quantized incremental algorithms for distributed optimization," \emph{IEEE J. Sel. Areas. Comm.}, vol. 23, no. 4, pp. 798--808,  2005.

\bibitem{CL10}
Y.~Cui, and V.K.N.~Lau, ``Convergence-optimal quantizer design of distributed contraction-based iterative algorithms with quantized message passing.", \emph{IEEE Trans. Signal Process.}, vol.~58, no.~10, pp.~5196--5205,  2010.

\bibitem{YXZR12}
D.~Yuan, S.~Xu, H.~Zhao and L.~Rong, ``Distributed dual averaging method for multi-agent optimization with quantized communication," \emph{Systems $\&$ Control Letters}, vol.~61, no.~11, pp.~1053-1061,  2012.

\bibitem{YH14}
P.~Yi and Y.~Hong, ``Quantized Subgradient Algorithm and Data-Rate Analysis for Distributed Optimization," \emph{IEEE Trans. Control Netw. Syst.}, vol.~1, no.~4, pp.~380-392,  2014.

\bibitem{NNA15}
E.~Nekouei, G.~N.~Nair and T.~Alpcan, ``Performance Analysis of Gradient-Based Nash Seeking Algorithms Under Quantization'', in press, \emph{IEEE Trans. Autom. Control}, 2016.


\bibitem{NNA15-CDC}
E.~Nekouei, G.~N.~Nair and T.~Alpcan, ``Convergence Analysis of Quantized Primal-dual Algorithm in Quadratic Network Utility Maximization Problems," 54th \emph{Int. Conf. on Decision and Control (CDC)}, pp. 2655--2660, 2015.

%\bibitem{CFSZ08}
%R.~Carli, F.~Fagnani, A.~Speranzon and S.~Zampieri, ``Communication constraints in the average consensus problem," \emph{Automatica}, vol.~44, no.~3, pp.~671-684, 2008.

%\bibitem{NFZE07}
%G.N.~Nair, F.~Fagnani, S.~Zampieri, and R.J.~Evans, ``Feedback Control Under Data Rate Constraints: An Overview," \emph{Proc. IEEE}, vol.~95, no.~1, pp.~108-137, Jan. 2007.

%\bibitem{Dong13}
%R.~Dong, and Z.~Geng, ``Design and analysis of quantizer for multi-agent systems with a limited rate of communication data," \emph{Communications in Nonlinear Science and Numerical Simulation}, pp. 282-290, vol. 18, no. 2, 2013.

%\bibitem{KBS07}
%A.~Kashyap, T.~Ba\c{s}ar and R.~Srikant, ``Quantized consensus," \emph{Automatica}, vol.~43, no.~7 pp.~1192-1203, 2007.

%\bibitem{NOOT09}
%A.~Nedi\'{c}, A.~Olshevsky, A.~Ozdaglar and J.~N.~Tsitsiklis, ``On distributed averaging algorithms and quantization effects," \emph{IEEE Trans. Autom. Control}, vol.~54, no.~11, pp.~2506-2517, Nov. 2009.

%\bibitem{YX11}
%K.~You and L.~Xie, ``Network Topology and Communication Data Rate for Consensusability of Discrete-Time Multi-Agent Systems," \emph{IEEE Trans. Autom. Control}, vol.~56, no.~10, pp.~2262-2275, Oct. 2011.

%\bibitem{MFDN09}
%P.~Minero, M.~Franceschetti, S.~Dey and G.N.~Nair, ``Data Rate Theorem for Stabilization Over Time-Varying Feedback Channels," \emph{IEEE Trans. Autom. Control}, vol.~54, no.~2, pp.243-255, Feb. 2009.

\bibitem{NE04}
G.~N.~Nair and R.J.~Evans, ``Stabilizability of stochastic linear systems with finite feedback data rates," \emph{SIAM J. Control Optim.}, vol.~43, no.~2, pp.~413--436, 2004.

\bibitem{FMS10}
J.S.~Freudenberg, R.H.~Middleton, V.~Solo, ``Stabilization and disturbance attenuation Over a Gaussian Communication Channel," \emph{IEEE Trans. Autom. Control}, vol.~55, no.~3, pp.795--799,  2010.

%\bibitem{NE04}
%G.N.~Nair and R.J.~Evans, ``Stabilizability of stochastic linear systems with finite feedback data rates," \emph{SIAM J. Control Optim.}, vol.~43, no.~2, pp.~413-436, 2004.

\bibitem{Cover}
T.~M.~Cover and J.~A.~Thomas, ``\emph{Elements of Information Theory}", 2nd ed., New York: Wiley, 2005.


\bibitem{leon-garcia}
A.~Leon-Garcia,``Probability and Random Processes for Electrical Engineering'', 2nd ed., Massachusetts: Addison-Wesley, 1994.

 \end{thebibliography}

\end{document}